\documentclass{amsart} 
\usepackage{amssymb,latexsym,textcomp}
\usepackage{xcolor}

\usepackage{cite} 

\usepackage{hyperref}
\hypersetup{colorlinks,allcolors=blue}
\usepackage[OT2,T1]{fontenc}

\theoremstyle{plain}
\newtheorem{theorem}{Theorem}[section]
\newtheorem{lemma}[theorem]{Lemma}
\newtheorem{corollary}[theorem]{Corollary}

\theoremstyle{definition}

\theoremstyle{remark}
\newtheorem{remark}{Remark}[section]

\DeclareMathOperator{\diam}{diam}
\DeclareMathOperator{\Ker}{Ker}

\newcommand{\C}{\mathbb C}

\newcommand{\strongperp}{\mathrel{\perp\!\!\!\!\perp}^S}
\newcommand{\strongp}{\mathrel{\perp}^S}
\newcommand\rk{\mathop{\mathrm{rk}}}

\def\FF{\mathbb F}

\def\Image{\mathop{\mathrm{Im}}}
\def\rank{\mathop{\mathrm{rk}}}
\def\Ker{\mathop{\mathrm{Ker}}}

\def\Span{\mathop{\mathrm{span}}}

\numberwithin{equation}{section} 

\begin{document}
\title[Eccentricity in mutual strong BJ orthographs of $C^\ast$-algebras]{Eccentricity in mutual strong Birkhoff-James orthographs of $C^\ast$-algebras} 

\author{Srdjan Stefanovi\'{c}}
\address{University of Belgrade\\ Faculty of Mathematics\\ Student\/ski trg 16-18\\ 11000 Belgrade\\ Serbia}

\email{srdjan.stefanovic@matf.bg.ac.rs}

\thanks{The research was supported by the Serbian Ministry of Education, Science and Technological Development through Faculty of Mathematics, University of Belgrade.}

\begin{abstract} We determine the eccentricity of every vertex in the orthograph of mutual strong Birkhoff-James orthogonality for arbitrary finite-dimensional $C^{\ast}$-algebras and for the algebra $K(H)$ of compact operators on an infinite-dimensional separable Hilbert space. This gives a more detailed description of the orthograph introduced by Aramba\v{s}i\'c et al. (Banach J. Math. Anal. 14 (2020), 1751-1772). As an application, we characterize scalar multiples of unitary elements among invertible elements of finite-dimensional $C^\ast$-algebras in terms of eccentricity alone. We also show that $M_3(\C)$ is the only finite-dimensional $C^\ast$-algebra whose orthograph has a non-trivial connected component of radius 3.
\end{abstract}




\maketitle

\section{Introduction and preliminaries}

Birkhoff-James orthogonality, originating in the work of Birkhoff \cite{Birkhoff} and further developed by James \cite{James1,James2, James3}, provides a notion of orthogonality in normed spaces that does not require an inner product. Let $X$ be a complex normed space. An element $x\in X$ is said to be Birkhoff-James orthogonal to $y\in X$, written $x\perp_{BJ}y$, if
$$\|x+\lambda y\|\geq\|x\|,\quad\text{for every }\lambda\in\mathbb C.$$
When the norm is induced by an inner product, this relation coincides with the usual orthogonality. In a general normed space, however, it need not be symmetric.

Graphs associated with orthogonality provide a useful framework for studying the geometry encoded by this relation. Symmetrized Birkhoff-James orthogonality in arbitrary normed spaces was studied by Aramba\v{s}i\'c et al.\ in \cite{Kuzma2021JMAA}. Guterman, Kuzma, Singla, and Zhilina \cite{BJNorm} showed that the directed Birkhoff-James orthogonality graph of a finite-dimensional normed space determines its dimension, smooth points and maximal faces of the unit ball. Related questions concerning the information contained
in Birkhoff-James orthogonality were investigated by Kuzma and Singla for simple, abelian, and general (real) finite-dimensional $C^\ast$-algebras \cite{KuzmaSushilFilomat,KuzmaSushilDebrecen,KuzmaSushilPAMS}. For a broader account of Birkhoff-James orthogonality, its characterizations, preservers and associated graphs, we refer to the survey by Aramba\v{s}i\'c, Guterman, Kuzma, and Zhilina \cite{ArambasicSurvey}.

A strong version of Birkhoff-James orthogonality was introduced by Aramba\v{s}i\'c and Raji\'c \cite{ArambasicAFA2014} in the setting of Hilbert $C^\ast$-modules. If $\mathcal M$ is a Hilbert module over a $C^\ast$-algebra $\mathfrak A$, then $a\in\mathcal M$ is strong Birkhoff-James orthogonal to $b\in\mathcal M$, written $a\strongp_{BJ}b$, if
$$\|a+bc\|\geq\|a\|,\quad\text{for every }c\in\mathfrak A.$$
Throughout this paper, we regard $\mathfrak A$ as a Hilbert $\mathfrak A$-module over itself. Since strong Birkhoff-James orthogonality need not be symmetric, we consider its mutual version:
$$a\strongperp_{BJ}b\quad\Longleftrightarrow\quad a\strongp_{BJ}b\ \text{and}\ b\strongp_{BJ}a.$$
The orthograph associated with this relation was introduced by Aramba\v{s}i\'c et al.\ in \cite{ArambasicBJMA2020}. For a $C^\ast$-algebra $\mathfrak A$, let $\Gamma(\mathfrak A)$ be the undirected graph whose vertex set is
$$V(\Gamma(\mathfrak A))=\{[a]=\mathbb Ca\mid a\in\mathfrak A\setminus\{0\}\}.$$
Two vertices $[a]$ and $[b]$ are adjacent if
$a\strongperp_{BJ}b$. This definition is independent of the choice of representatives, since mutual strong Birkhoff-James orthogonality is invariant under multiplication by non-zero scalars. We shall therefore usually write $a$ instead of $[a]$.

In the paper \cite{ArambasicBJMA2020}, the authors studied isolated vertices and the diameter of the orthograph for $B(H)$ - the algebra of bounded operators on a complex Hilbert space, as well as for commutative $C^\ast$-algebras. In \cite{Keckic2023JMAA}, Kečkić and Stefanović answered the questions posed in the preceding paper and characterized isolated vertices in arbitrary, possibly non-unital, $C^\ast$-algebras. They also proved that, except for the three small algebras $\mathbb C$, $\mathbb C\oplus\mathbb C$, and $M_2(\mathbb C)$, there is a unique non-trivial connected component, and obtained bounds on its diameter. Finally, in \cite{Stefanovic}, the diameter of an arbitrary finite-dimensional $C^\ast$-algebra was determined. We will repeatedly use the main result of that paper: apart from the three exceptional algebras mentioned above, the unique non-trivial connected component has diameter $4$ precisely for $M_3(\mathbb C)$ and $\mathbb C\oplus M_2(\mathbb C)$, and diameter $3$ in every other finite-dimensional case.

For convenience, we recall the graph-theoretic notions used throughout the paper. Unless otherwise specified, all graph-theoretic quantities are computed within the connected component of the orthograph. For vertices $x,y$ in a connected graph $\Gamma$, the distance $d(x,y)$ is the length of a shortest path joining them. The eccentricity of $x$ is
$$e(x)=\sup\{d(x,y)\mid y\in V(\Gamma)\}.$$
The diameter and radius are then defined, respectively, as
$$\diam(\Gamma)=\sup\{e(x)\mid x\in V(\Gamma)\} \quad\text{and}\quad r(\Gamma)=\inf\{e(x)\mid x\in V(\Gamma)\}.$$
The vertices with the lowest eccentricity form the center of the graph. 

Graph-theoretic methods have also been used to characterize maps preserving the corresponding relations. We mention two results that provide a starting point for our study. Dolinar and Kuzma \cite{DolinarKuzma2016} used vertices at maximal distance in the commuting graph to describe surjective maps on $M_n(\FF), n\geq5$ that preserve commutativity in one direction (without assuming additivity or linearity) and map non-scalar matrices to non-scalar matrices; here $\FF$ is an algebraically closed field. The second paper we would like to highlight, especially since it studies the same relation as we do, is that of Aramba\v{s}i\'c et al.\ \cite{StrongMutualPreservers}, which is a natural continuation of the research begun in \cite{ArambasicBJMA2020}. In this paper, the authors classified all surjective linear maps on $B(H)$ preserving mutual strong Birkhoff-James orthogonality. 

The aim of this paper is to describe more detailed the orthograph of mutual strong Birkhoff-James orthogonality and to improve some of the results obtained in \cite{ArambasicBJMA2020,Stefanovic}. More precisely, we determine the eccentricity of every vertex in the orthograph of an arbitrary finite dimensional $C^\ast$-algebra. For compact operators, we prove that the non-trivial component has diameter $3$ and obtain a complete eccentricity classification. The descriptions relate eccentricity to norm-attaining subspaces and kernels of adjoints.

The paper is organized as follows. In Section 2, we treat $M_n(\mathbb C)$ for $n\geq3$, distinguishing the exceptional case $n=3$ from $n\geq4$. We treat the simple finite-dimensional case separately for two reasons. It illustrates the main ideas used in the general case, and the results for $M_n(\C), n\geq4$, extend naturally to $K(H)$, where $H$ is an infinite-dimensional separable Hilbert space. Section 3 deals with finite-dimensional $C^\ast$-algebras having at least two summands, beginning with $\mathbb C\oplus M_2(\mathbb C)$ and then establishing the general classification. We also give examples and describe the radius and center as a consequence. Among the finite-dimensional $C^\ast$-algebras not isomorphic to any of the three small algebras mentioned above, the radius is $3$ precisely for $M_3(\mathbb C)$ and $2$ otherwise. We also characterize scalar multiples of unitary elements among all invertible elements of a finite-dimensional $C^\ast$-algebra. In Section 4, we study the orthograph of $K(H)$, where $H$ is an infinite-dimensional separable Hilbert space. We determine the eccentricity of every vertex and hence the radius and diameter of the unique non-trivial connected component.

\section{The case of \texorpdfstring{$M_n(\mathbb C)$, $n\geq3$}
{Mn(C), n >= 3}}

We distinguish between two cases $n=3$ and  $n\geq4$, since the corresponding non-trivial connected components have different diameters. In both cases, for a given $A\in M_n(\C)$, the two subspaces that will play a central role are 
$$AM_0(A)\ \text{and}\ \Ker A^*=(\Image A)^{\perp}.$$ 
Here,
$$M_0(A):=\Ker (A^\ast A-\|A\|^2 I)=\{x\in\C^n\mid \|Ax\|=\|A\|\|x\|\}$$
denotes the norm-attaining subspace of $A$.

We will repeatedly use the following characterization, known in the literature as the strong version of the Stampfli, Magajna, Bhatia-\v{S}emrl theorem (see \cite{Stampfli}, \cite{Magajna1993} and \cite{BhatiaSemrl}). We will use a version of this theorem from \cite[Example 4.1]{RACSAM6}. We note that we also include the case $n=1$, which will be needed in the next section. This case is immediate, since for $A,B\in\C$ with $A\neq0$, both conditions in the theorem hold if and only if $B=0$.

\begin{theorem}\label{glavna}
    Let $A\in M_n(\C)\setminus\{0\}$ and $B\in M_n(\C), n\geq 1$. Then $A\strongp_{BJ}B$ if and only if 
    $$AM_0(A)\cap(\Image B)^{\perp}\neq\{0\}.$$
\end{theorem}

We begin by proving that there is no vertex with eccentricity equal to 1.

\begin{lemma}\label{e>1}
    Let $A\in M_n(\C),n\geq3$ be a non-zero non-invertible matrix. Then, there exists a matrix $B\in M_n(\C)$ such that $ A\not\!\strongperp_{BJ}B$. In particular, $e(A)\geq2$.
\end{lemma}
\begin{proof}
    Assume $\|A\|=1$ and $$A=x_1y_1^*+\sum_{i=2}^{n-1}\sigma_ix_iy_i^*,\quad 1\geq\sigma_2\geq\dots\geq\sigma_{n-1}\geq0$$
    is its SVD. Define $$B:=x_1y_2^*.$$ Then $[B]\neq[A]$ and $B$ attains its norm only on scalar multiples of $y_2$ and $A^*By_2=A^*x_1=y_1\neq0$, so by Theorem \ref{glavna} $B\not\!\strongp_{BJ}A$. Also, $B$ is a non-zero, non-invertible matrix, so it belongs to the same unique non-trivial component as $A$, hence $e(A)\geq2$.
\end{proof}

\subsection{The case of \texorpdfstring{$M_3(\mathbb C)$}{M3(C)}}

First, we will show that in the case of $M_3(\C)$ there are no vertices of eccentricity 2.

\begin{lemma}\label{radius M_3}
 The radius of the non-trivial component of $\Gamma(M_3(\C))$ is greater than or equal to 3.
\end{lemma}

\begin{proof}
    Fix a non-zero non-invertible $A\in M_3(\C)$ and assume $\|A\|=1$. We consider two cases. 

    \underline{\emph{Case 1: $\rank A=1$.}} Then $A=x_1y_1^\ast$ is its SVD. Extend $x_1$ to an orthonormal basis $\{x_1,x_2,x_3\}$ of $\mathbb{C}^3$. Define $$B:=\frac{x_1+x_2}{\sqrt{2}}\left(\frac{x_1+x_2}{\sqrt{2}}\right)^\ast+x_3x_3^{\ast}.$$
    Then $(\Image B)^{\perp}=\Span\left\{\frac{x_1-x_2}{\sqrt{2}}\right\}$, and $A$ attains its norm only on scalar multiples of $y_1$, so $AM_0(A)=\Span\{x_1\}$. It follows that $AM_0(A)\cap (\Image B)^{\perp}=\{0\}$, so $A\not\!\strongp_{BJ} B$.

    Suppose, toward a contradiction, that there exists non-zero $C\in M_3(\C)$ such that $A\strongperp_{BJ}C\strongperp_{BJ}B$. By Theorem \ref{glavna}, we conclude $$\Span\{x_1\}\cap\left(\Image C\right)^{\perp}\neq\{0\}\ \text{and}\ CM_0(C)\cap\Span\left\{\frac{x_1-x_2}{\sqrt{2}}\right\}\neq\{0\}.$$
    Then $\frac{x_1-x_2}{\sqrt{2}}\in CM_0(C)\subseteq\Image C$ and $x_1\in (\Image C)^{\perp}$. This is impossible since $\langle x_1,\frac{x_1-x_2}{\sqrt{2}}\rangle=\frac{1}{\sqrt{2}}\neq0$. Consequently, $d(A,B)\geq3$. 

    \underline{\emph{Case 2: $\rank A=2$.}} Then $A=x_1y_1^*+\sigma x_2y_2^{*}$, where $0<\sigma\leq1$. Then $\left(\Image A\right)^{\perp}=\Span\{x_3\}$. Define $$B:=\frac{x_3+x_1}{\sqrt{2}}\left(\frac{x_3+x_1}{\sqrt{2}}\right)^*.$$
    Then $BM_0(B)\cap(\Image A)^{\perp}=\{0\}$, so by Theorem \ref{glavna} $B\not\!\strongp_{BJ}A$. Suppose again that there exists a non-zero $C$ that is mutual strong BJ orthogonal to both $A$ and $B$. By the same theorem and the fact that $BM_0(B)=\Span{\left\{\frac{x_3+x_1}{\sqrt{2}}\right\}}$ we obtain
    $$x_3\in CM_0(C)\subseteq\Image C\ \text{and}\ \frac{x_3+x_1}{\sqrt{2}}\in(\Image C)^{\perp},$$
    which is impossible.
\end{proof}

Now, from the previous lemma, Lemma \ref{e>1} and the fact that diameter of $M_3(\C)$ is equal to 4, we know that all non-isolated vertices have eccentricity 3 or 4. Our next goal is to classify the eccentricity of vertices. In the next two lemmas, we will describe the vertices of eccentricity 3.

\begin{lemma}
    If $A\in M_3(\C)$ is a rank-one operator, then $e(A)=3$.
\end{lemma}
\begin{proof}
    Lemma \ref{radius M_3} yields $e(A)\geq3$. After normalization, we may assume that $\|A\|=1$. Let 
    $$A=x_1y_1^*$$
    be its SVD. Take any non-zero non-invertible $B\in M_3(\C)$, assume also that $\|B\|=1$ and let 
    $$B=t_1z_1^*+\sigma_2t_2z_2^*,\quad1\geq\sigma_2\geq0,$$
    be its SVD. If $\{t_1,t_2,t_3\}$ is an orthonormal basis of $\C^3$, define $$C_2:=t_3t_3^*.$$ Then $C_2\strongperp_{BJ}B$. Namely, $C_2$ attains its norm at the vector $t_3$ and $B^*C_2t_3=B^*t_3=0$; on the other hand, $B$ attains its norm at the vector $z_1$ and $C_2^*Bz_1=C_2t_1=0$. Now, let $v\in\Span\{x_1,t_3\}^{\perp}$ be a unit vector and  define $C_1$ as 
    $$C_1:=vv^*.$$
    Then it holds $A\strongperp_{BJ}C_1\strongperp_{BJ}C_2\strongperp_{BJ}B,$
    so $e(A)=3$. It remains only to verify the orthogonalities involving $C_1$. $A$ attains its norm at the vector $y_1$, and $C_1^*Ay_1=C_1x_1=0$, so $A\strongp_{BJ} C_1$. Also, $C_1$ attains its norm at vector $v$ so $A^*C_1v=A^*v=\langle x_1,v\rangle y_1=0$, so $C_1\strongp_{BJ}A$. Finally, the fact that $\langle v,t_3\rangle=0$ gives us $C_1\strongperp_{BJ} C_2$.
\end{proof}

\begin{lemma}
    If $A\in M_3(\C)$ is a rank-two operator satisfying $AM_0(A)=\Image A$, then $e(A)=3$.
\end{lemma}
\begin{proof}
    Lemma \ref{radius M_3} yields $e(A)\geq3$. After normalization, we may assume that $\|A\|=1$. Let 
    $$A=x_1y_1^*+x_2y_2^*$$
    be its SVD. Take any non-zero non-invertible $B\in M_3(\C)$, assume also that $\|B\|=1$ and let 
    $$B=t_1z_1^*+\sigma_2t_2z_2^*,\quad1\geq\sigma_2\geq0,$$
    be its SVD. If $\{x_1,x_2,x_3\}$ and $\{t_1,t_2,t_3\}$ are orthonormal bases of $\C^3$, then we can define $$C_1:=P_{\Span\{x_3,t_1\}}$$ as the orthogonal projection onto the subspace $\Span\{x_3,t_1\}$ and $$C_2:=t_3t_3^*+t_2t_2^*,$$ then $A\strongperp_{BJ}C_1\strongperp_{BJ}C_2\strongperp_{BJ}B.$ Namely, $$AM_0(A)\cap(\Image C_1)^{\perp}=\Span\{x_1,x_2\}\cap x_3^{\perp}\cap t_1^{\perp}=\Span\{x_1,x_2\}\cap t_1^{\perp}\neq\{0\},$$ 
    by a dimension argument. Also, $C_1\strongp_{BJ}A$ by the fact that $C_1$ attains its norm at $x_3$ and  $A^*C_1x_3=0$. Further, $C_1$ attains its norm at $t_1$ and $C_2^*C_1t_1=0$ so $C_1\strongp_{BJ}C_2$. Again,  $$C_2M_0(C_2)\cap(\Image C_1)^{\perp}=\Span\{t_3,t_2\}\cap t_1^{\perp}\cap x_3^{\perp}=\Span\{t_3,t_2\}\cap x_3^{\perp}\neq\{0\},$$  so $C_2\strongp_{BJ}C_1$. Finally, since $B$ attains its norm at $z_1$ and $C_2^*Bz_1=C_2^*t_1=0$ we get $B\strongp_{BJ}C_2$; also $C_2$ attains its norm at $t_3$ and $B^*C_2t_3=B^*t_3=0$, so $C_2\strongp_{BJ}B$. This completes the proof.
\end{proof}

Finally, we characterize the vertices whose eccentricity equals the diameter of connected component of $\Gamma(M_3(\C))$, that is, the vertices having another vertex at distance 4.

\begin{lemma}
    Let $A\in M_3(\C)$ be a rank-two operator such that $\dim (AM_0(A))=1$. Then $e(A)=4$.
\end{lemma}
\begin{proof}
    From Lemma \ref{radius M_3}, we know that $e(A)\geq3$. Assume that $\|A\|=1$ and let 
    $$A=x_1y_1^*+\sigma_2x_2y_2^*,\quad 1>\sigma_2>0$$
    be its SVD and let $\{x_1,x_2,x_3\}$ be an orthonormal basis of $\C^3$. Define 
    $$B=\left(\frac{x_1+x_3}{\sqrt{2}}\right)\left(\frac{x_1+x_3}{\sqrt{2}}\right)^*+\frac{1}{2}x_2x_2^*.$$
    We claim that $d(A,B)=4$. 

    First, $A\not\!\strongp_{BJ}B$ by Theorem \ref{glavna} because $A$ attains its norm only on scalar multiples of $y_1$ and $B^*Ay_1=B^*x_1=\frac{x_1+x_3}{2}\neq0$.

    Next, we will prove that the distance is not 2. Suppose, to the contrary, that there exists non-zero $C\in M_3(\C)$ such that $A\strongperp_{BJ}C\strongperp_{BJ}B$. From $A\strongp_{BJ}C$ and Theorem \ref{glavna} we can conclude that $x_1\in(\Image C)^{\perp}$. Also, by the same theorem and $C\strongp_{BJ}B$ we have $\frac{x_1-x_3}{\sqrt{2}}\in CM_0(C)\subseteq\Image C$, which is not possible because $\langle x_1,\frac{x_1-x_3}{\sqrt{2}}\rangle=\frac{1}{\sqrt{2}}\neq0$.
     
    Finally, suppose that there exist non-zero $C_1,C_2\in M_3(\C)$ such that $A\strongperp_{BJ}C_1\strongperp_{BJ}C_2\strongperp_{BJ}B$. Then by Theorem \ref{glavna} and facts that $\Image A=\Span\{x_1,x_2\}$, $AM_0(A)=\Span\{x_1\}$, $\Image B=\Span\left\{\frac{x_1+x_3}{\sqrt{2}},x_2\right\}$ and $BM_0(B)=\Span\left\{\frac{x_1+x_3}{\sqrt{2}}\right\}$ we conclude
    $$x_3\in C_1M_0(C_1)\ , \ x_1\in(\Image C_1)^{\perp}\ ,\ \frac{x_1-x_3}{\sqrt{2}}\in C_2M_0(C_2)\ \text{and}\ \frac{x_1+x_3}{\sqrt{2}}\in(\Image C_2)^{\perp}.$$
    Consequently, after normalizing $C_1$ and $C_2$, their SVDs have the form 
    $$C_1=x_3t_3^*+\lambda_2x_2t_2^*\ \text{and}\ C_2=\frac{x_1-x_3}{\sqrt{2}}z_3^*+\tau_2x_2z_2^*,\quad 0\leq\lambda_2,\tau_2\leq1,$$
    for some orthonormal pairs of vectors in $\C^3$: $\{t_2,t_3\}$ and $\{z_2,z_3\}$. If $\lambda_2<1$, then $C_1$ attains its norm only on scalar multiples of $t_3$ and $$C_2^*C_1t_3=C_2^*x_3=-\frac{1}{\sqrt{2}}z_3\neq0,$$
    so $C_1\not\!\strongp_{BJ}C_2$, thus $\lambda_2=1$. A symmetric argument shows that $\tau_2=1$. So $C_1$ attains its norm on every vector $v\in\Span\{t_3,t_2\}$. But, if $v=\alpha t_3+\beta t_2$ for some $\alpha,\beta\in\C$ such that $|\alpha|^2+|\beta|^2=1$, we get 
    $$C_2^*C_1v=\left(z_3\left(\frac{x_1-x_3}{\sqrt{2}}\right)^*+z_2x_2^*\right)(\alpha x_3+\beta x_2)=-\frac{\alpha}{\sqrt{2}}z_3+\beta z_2\neq0.$$
    This contradicts $C_1\strongp_{BJ} C_2$. Therefore, no path of length at most 3 joins $A$ and $B$. Since the diameter is 4, we conclude that $e(A)=4$.
\end{proof}

We have therefore proved the following theorem.

\begin{theorem}\label{M_3}
    The unique non-trivial connected component of $\Gamma(M_3(\C))$ has radius 3 and diameter 4. Its center consists precisely of the projective classes of non-zero singular matrices $A\in M_3(\C)$ such that either $\rk A=1$, or $\rk A=2$ and $AM_0(A)=\Image A$. 
\end{theorem}

\subsection{The case of \texorpdfstring{$M_n(\mathbb C)$, $n\geq4$}
{Mn(C), n >= 4}}

In the case when $n\geq4$, we know that the diameter of  $\Gamma(M_n(\C))$ is equal to 3, and by Lemma \ref{e>1} we know that the eccentricity of each non-isolated vertex is at least 2. Therefore, a non-isolated vertex in $M_n(\C)$ can have only two possible eccentricities: 2 or 3. First, we will show one class of vertices with eccentricity 3.

\begin{lemma}\label{AM_0(A) jednodimenzionalan}
    If $A\in M_n(\C),n\geq4$ is non-invertible element and $\dim (AM_0(A))=1$, then $e(A)=3$.
\end{lemma}
\begin{proof}
    Assume that $\|A\|=1$ and let 
    $$A=x_1y_1^*+\sum_{i=2}^{n-1}\sigma_ix_iy_i^*,\quad 1>\sigma_2\geq\dots\geq\sigma_{n-1}\geq0$$
    be its SVD. Let $\{x_1,\dots,x_{n-1},x_n\}$ be an orthonormal basis for $\C^n$. Define 
    \begin{equation}
    B=\frac{x_n+x_1}{\sqrt{2}}\left(\frac{x_n+x_1}{\sqrt{2}}\right)^*+\frac{1}{2}\sum_{i=2}^{n-1}x_ix_i^*.
    \end{equation}
    We will prove that $d(A,B)\geq3$, thus we conclude that $e(A)=3$, since the diameter of the non-trivial component is 3.

    First, $A$ attains its norm only on scalar multiples of $y_1$ and $B^*Ay_1=B^*x_1=\frac{x_1+x_n}{2}\neq0$, so $A\not\!\strongp_{BJ} B$.

    If there exists $C\in M_n(\C)\setminus\{0\}$ such that $A\strongperp_{BJ}C\strongperp_{BJ}B$, then by Theorem \ref{glavna}, $AM_0(A)\cap(\Image C)^{\perp}\neq\{0\}$. Since $AM_0(A)=\Span\{x_1\}$, it follows that $x_1\in(\Image C)^{\perp}$. But by the same theorem, $CM_0(C)\cap(\Image B)^{\perp}\neq\{0\}$ and applying $(\Image B)^{\perp}=\Span\left\{{\frac{x_1-x_n}{\sqrt{2}}}\right\}$, we get $\frac{x_1-x_n}{\sqrt{2}}\in CM_0(C)\subseteq\Image C$, which is impossible because $\left\langle x_1,\frac{x_1-x_n}{\sqrt{2}}\right\rangle=\frac{1}{\sqrt{2}}\neq0$.
\end{proof}

\begin{remark}\label{remark:one-sided-AM0}
    The proof of Lemma \ref{AM_0(A) jednodimenzionalan} yields a stronger, one-sided statement. More precisely, let $n\geq2$ and let $A\in M_n(\C)\setminus\{0\}$, after normalization, has SVD of the form
    $$A=x_1y_1^*+\sum_{i=2}^{n}\sigma_ix_iy_i^*\ 1>\sigma_2\geq\dots\geq\sigma_n\geq0,$$
    so $AM_0(A)=\Span\{x_1\}$,
    and define
    $$B=\frac{x_n+x_1}{\sqrt2}
        \left(\frac{x_n+x_1}{\sqrt2}\right)^*
        +\frac12\sum_{i=2}^{n-1}x_ix_i^*.$$
    Then $A\not\!\strongp_{BJ}B$ and there is no non-zero $C\in M_n(\C)$ such that $A\strongp_{BJ}C\strongp_{BJ}B.$ Indeed,
    $(\Image B)^\perp=\Span\left\{\frac{x_1-x_n}{\sqrt2}\right\},$ so $AM_0(A)\cap(\Image B)^\perp=\{0\}$. Furthermore, if $A\strongp_{BJ}C$, then
            $$x_1\in(\Image C)^\perp,$$
    whereas $C\strongp_{BJ}B$ implies
        $$\frac{x_1-x_n}{\sqrt2}\in CM_0(C)\subseteq\Image C.$$
    This is impossible because $\left\langle x_1,\frac{x_1-x_n}{\sqrt2}\right\rangle
        =\frac1{\sqrt2}\neq0.$
    
    Thus, the proof does not use either
    $C\strongp_{BJ}A$ or $B\strongp_{BJ}C$. For $n=2$, the sum in the definition of $B$ is understood to be empty. Hence the statement applies, in particular, to both $n=2$ and $n=3$. Also, here we allow $A$ to be invertible, which will be important for us in further application.
\end{remark}

There is only one more family of vertices of eccentricity 3.

\begin{lemma}\label{ImA n-1-dimenzionalan}
    If $A\in M_n(\C),n\geq4$ and $\dim (\Image A)^{\perp}=1$, then $e(A)=3$.    
\end{lemma}
\begin{proof}
    The proof is similar to that of Lemma \ref{AM_0(A) jednodimenzionalan}. Suppose  $\|A\|=1$ and let $$A=x_1y_1^*+\sum_{i=2}^{n-1}\sigma_ix_iy_i^*,\ 1\geq\sigma_2\geq\dots\geq\sigma_{n-1}>0$$
    be its SVD. Let $\{x_1,\dots,x_{n-1},x_n\}$ be an orthonormal basis for $\C^n$. Define
    $$B=\frac{x_n+x_1}{\sqrt2}\left(\frac{x_n+x_1}{\sqrt2}\right)^*.$$
    Then $B$ attains its norm only on scalar multiples of $\frac{x_n+x_1}{\sqrt{2}}$ and $A^*B\left(\frac{x_n+x_1}{\sqrt{2}}\right)=A^*\left(\frac{x_n+x_1}{\sqrt{2}}\right)=\frac{1}{\sqrt{2}}y_1\neq0$, so $B\not\!\strongp_{BJ}A$.

    Again, if there exists $C\in M_n(\C)\setminus\{0\}$ such that $A\strongperp_{BJ}C\strongperp_{BJ}B$, then by Theorem \ref{glavna} we deduce $CM_0(C)\cap(\Image A)^{\perp}\neq\{0\}$. By the hypothesis $(\Image A)^{\perp}=\Span\{x_n\}$, so $x_n\in CM_0(C)\subseteq\Image C$. By the same theorem, we also know that $(\Image C)^{\perp}\cap BM_0(B)\neq\{0\}$, and by definition of $B$, $BM_0(B)=\Span\{\frac{x_1+x_n}{\sqrt{2}}\}$, so $\frac{x_1+x_n}{\sqrt{2}}\in(\Image C)^{\perp}$. But this is a contradiction because $\left\langle x_n,\frac{x_1+x_n}{\sqrt{2}}\right\rangle=\frac{1}{\sqrt{2}}\neq0$. Therefore, no path of length at most 2 joins $A$ and $B$. Since the diameter is 3, we conclude that $e(A)=3$.
\end{proof}

\begin{remark}\label{remark:one-sided-kernel}
    The proof of Lemma \ref{ImA n-1-dimenzionalan} also yields a stronger, one-sided statement. Let $n\geq2$ and let $A\in M_n(\C)$, after normalization, has SVD of the form
    $$A=\sum_{i=1}^{n-1}\sigma_ix_iy_i^*,\ \sigma_1\geq\dots\geq\sigma_{n-1}>0,$$
    so $\dim(\Image A)^\perp=1.$ Let $\{x_1,\dots,x_{n-1}, x_n\}$ be an orthonormal basis for $\C^n$. Then 
    $(\Image A)^\perp=\Span\{x_n\}$. Define
    $$B=\frac{x_n+x_1}{\sqrt2}\left(\frac{x_n+x_1}{\sqrt2}\right)^*.$$
    
    Then $B\not\!\strongp_{BJ}A$ and there is no non-zero $C\in M_n(\C)$ such that $B\strongp_{BJ}C\strongp_{BJ}A.$ Indeed,
    $BM_0(B)=\Span\left\{\frac{x_1+x_n}{\sqrt2}\right\},$
    and hence $BM_0(B)\cap(\Image A)^\perp=\{0\}.$ Moreover, $B\strongp_{BJ}C$ implies $\frac{x_1+x_n}{\sqrt2}\in(\Image C)^\perp,$
    while $C\strongp_{BJ}A$ implies
    $x_n\in CM_0(C)\subseteq\Image C.$ This is impossible because $\left\langle x_n,\frac{x_1+x_n}{\sqrt2}\right\rangle=\frac1{\sqrt2}\neq0.$
    Thus, the proof does not use either $A\strongp_{BJ}C$ or $C\strongp_{BJ}B$. Again the statement covers both $n=2$ and $n=3$.
\end{remark}

The following lemma shows that all remaining vertices have eccentricity 2. 

\begin{lemma}\label{oba bar dva}
    If $A\in M_n(\C), n\geq4$ with $\dim AM_0(A)\geq2$ and $\dim (\Image A)^{\perp}\geq2$, then $e(A)=2$.
\end{lemma}
\begin{proof}
    We can assume that $\|A\|=1$ and let $$A=\sum_{i=1}^{n-2}\sigma_ix_iy_i^*,\quad 1=\sigma_1=\sigma_2\geq\dots\geq\sigma_{n-2}\geq0$$
    be its SVD. Take an arbitrary non-zero non-invertible $B$ (without loss of generality $\|B\|=1$) with SVD
    $$B=t_1z_1^*+\sum_{i=2}^{n-1}\tau_it_iz_i^*,\quad 1\geq\tau_2\geq\dots\geq\tau_{n-1}\geq0.$$
    Let $\{t_1,\dots,t_{n-1},t_n\}$ be an orthonormal basis for $\C^n$. Furthermore, let $v'$ be a unit vector in $\Span\{y_1,y_2\}$ and $v=Av'\in AM_0(A)$ be such that $v\perp t_n$. Such a vector $v$ exists by a dimension argument. Analogously, we can find a unit vector $w$ such that $w\in(\Image A)^{\perp}\cap t_1^{\perp}$ because $\dim(\Image A)^\perp\geq2$ and $\dim(t_1^\perp)=n-1$. Also, we have $\langle v,w\rangle=0$. Now, we can define $C$ as an orthogonal projection onto the subspace $\Span\{t_n,w\}$. We will prove that $A\strongperp_{BJ} C\strongperp_{BJ} B$.
    
    The operator $A$ attains the norm on $v'$ and $C^*Av'=Cv=0$ because $v\perp t_n$ and $v\perp w$, so $A\strongp_{BJ}C$. In the other direction, $C$ attains its norm at $w$ and $A^*Cw=A^*w=0$ because $w\in(\Image A)^{\perp}$, so $C\strongp_{BJ}A$.

    Furthermore, $B$ attains its norm at $z_1$ and $C^*Bz_1=Ct_1=0$ because $t_1\perp w$ and $t_1\perp t_n$, so $B\strongp_{BJ}C$. On the other hand, $C$ attains its norm at $t_n$ and $B^*Ct_n=B^*t_n=0$ because $t_n\in(\Image B)^\perp$, so we can finally conclude $C\strongp_{BJ}B$. Thus, every non-zero non-invertible $B$ is at distance at most 2 from $A$, so by Lemma \ref{e>1}, we have proven $e(A)=2$.
\end{proof}

We have therefore proved the following theorem.

\begin{theorem}\label{M_n,n>3}
    The unique non-trivial connected component of $\Gamma(M_n(\C)), n\geq4$ has radius 2 and diameter 3. Its center consists precisely of the projective classes of non-zero singular matrices $A\in M_n(\C)$ satisfying $$\dim(AM_0(A))\geq2\quad\text{and}\quad\dim(\Image A)^{\perp}\geq2.$$ 
\end{theorem}

\section{Cases with two or more summands}

First, we will prove that even in the case of two or more summands, there is no vertex of eccentricity 1.

\begin{lemma}\label{ecc>1 for all}
    Let $\mathfrak{A}=\bigoplus_{i=1}^{k}M_{n_i}(\C), k\geq2$ be a finite-dimensional $C^\ast$-algebra not isomorphic to $\C\oplus\C$. Then for every non-zero, non-invertible $a\in\mathfrak{A}$ its eccentricity is at least 2.
\end{lemma}
\begin{proof}
    Recall that, for every such finite-dimensional $C^\ast$-algebra, all non-zero non-invertible elements belong to the same non-trivial connected component. Assume $a=(A_1,\dots,A_k)$ and $\|a\|=1$. Suppose first that there exists an index $j$ such that $n_j\geq2$ and $A_j\in M_{n_j}(\C),\|A_j\|=1$. If $$A_j=x_{1}y_1^*+\sum_{i=2}^{n_j}\sigma_ix_iy_i^*,\quad 1\geq\sigma_2\geq\dots\geq\sigma_{n_j}\geq0,$$
    is SVD of $A_j$, define $b=(B_1,\dots,B_k)$ such that $B_j=x_1y_2^*$ and $B_i=0$ for $i\neq j$. Then, $B_jM_0(B_j)=\Span\{x_1\}\subseteq\Image A_j$, so $B_jM_0(B_j)\cap(\Image A_j)^{\perp}=\{0\}$. According to Theorem \ref{glavna}, we have $B_j\not\!\strongp_{BJ}A_j$, and from that $b\not\!\strongp_{BJ}a$. In the remaining case, $a$ attains its norm only on summands $\C$. After rescaling and reordering the summands, we can assume that 
    $$a=(1,a'), a'\in\mathfrak{B},$$
    where $\mathfrak{B}$ is the corresponding finite-dimensional $C^\ast$-algebra different from $\C$. Then $(1,0)\not\!\strongp_{BJ}a$ because 
    $$\|(1,0)+(1,a')(-1,0)\|=0<1=\|(1,0)\|.$$
    Therefore, we have shown that $e(a)\geq2$ but under the condition that $(1,a')$ is not a scalar multiple of $(1,0)$. So, the only remaining case is that $a=(1,0)$. But we actually solved this case too, because according to the previous $(1,0)\not\!\strongp_{BJ}(1,c)$ for any non-zero, non-invertible $c\in\mathfrak{B}$ (such $c$ exists because $\mathfrak{B}$ is different from $\C$).
\end{proof}

Next, we will prove the following very useful lemma.

\begin{lemma}\label{prelazak na istu normu}
    Let $\mathfrak{A}$ and $\mathfrak{B}$ be unital $C^\ast$-algebras, and let $a\in\mathfrak{A}$ and $b\in\mathfrak{B}$ satisfy $\|b\|>\|a\|>0$. If $b$ is not right invertible element, then there exists $b'\in\mathfrak{B}$ such that $b'\strongperp_{BJ}b$, $\|b'\|=\|a\|$ and
    $$(a,b)\strongperp_{BJ}(a,b').$$
\end{lemma}
\begin{proof}
    From the fact that $b$ is not right invertible, we know from \cite[Proposition 2.4]{Keckic2023JMAA} that there exists $b'$ such that 
    $$b'\strongperp_{BJ}b\quad\text{and}\quad\|b'\|=\|a\|.$$
    Then for every $(x,y)\in\mathfrak{A}\oplus\mathfrak{B}$ it holds
    $$\max\{\|a+ax\|,\|b+b'y\|\}\geq\|b+b'y\|\geq\|b\|=\max\{\|a\|,\|b\|\}$$
    and
    $$\max\{\|a+ax\|,\|b'+by\|\}\geq\|b'+by\|\geq\|b'\|=\max\{\|a\|,\|b'\|\},$$
    thus the desired conclusion follows.
\end{proof}

\subsection{Case $\C\oplus M_2(\C)$}

In this case we know the diameter is 4, and by Lemma \ref{ecc>1 for all} we know that the eccentricity of each non-isolated vertex is at least 2, therefore, for a given non-isolated vertex, we have that its possible eccentricity is 2, 3 or 4. Also, we have four types of non-invertible vertices: 
$$(1,0);\quad(0,A),\rk A=1;\quad(0,A),\rk A=2;\quad (1,A),\rk A=1.$$

Let us denote by $$A^{\strongp_{BJ}}=\{B\in\mathfrak{A}\mid A\strongp_{BJ}B\}$$ the outgoing neighborhood of $A$.

The following easy characterization will be useful for testing orthogonality in this exceptional case $\C\oplus M_2(\C)$.

\begin{lemma}\label{C+M_2}
    Let $(\lambda,A)\in\C\oplus M_2(\C)$ be a non-zero non-invertible element. Then we have three cases:
    \begin{itemize}
        \item[i)] If $|\lambda|<\|A\|$, then $(\lambda,A)^{\strongp_{BJ}}=\{(\mu,B)\in\C\oplus M_2(\C)\mid A\strongp_{BJ}B\}$.
        \item[ii)] If $|\lambda|>\|A\|$, then $(\lambda,A)^{\strongp_{BJ}}=\{(0,B)\in\C\oplus M_2(\C)\mid B \text{ arbitrary}\}$.
        \item[iii)] If $|\lambda|=\|A\|$, then $(\lambda,A)^{\strongp_{BJ}}=\{(\mu,B)\in\C\oplus M_2(\C)\mid A\strongp_{BJ}B\ \text{or}\ \mu=0\}$.
    \end{itemize}
\end{lemma}
\begin{proof}
    If $(\lambda,A)\strongp_{BJ}(\mu,B)$ by definition we have 
    $$\max\{|\lambda+\mu x|,\|A+BY\|\}\geq\max\{|\lambda|,\|A\|\},\ \text{ for all }(x,Y)\in\C\oplus M_2(\C).$$
    If $|\lambda|<\|A\|$, by setting $x=0$, we get $A\strongp_{BJ}B$. If $|\lambda|>\|A\|$ holds and $\mu\neq0$, we can put $x=-\mu^{-1}\lambda$ and $Y=0$ to get a contradiction, hence $\mu=0$. Finally, suppose that $|\lambda|=\|A\|$ holds. Assume the contrary, $A\not\!\strongp_{BJ}B$ and $\mu\neq0$. Then there exists $Y$ such that $\|A+BY\|<\|A\|$ and we can choose $x=-\mu^{-1}\lambda$ to get a contradiction. The converse implications in all three cases are trivial.  
\end{proof}

First, we will find an example of a vertex of eccentricity 2.

\begin{lemma}
    The eccentricity of $(1,0)$ is equal to 2.
\end{lemma}
\begin{proof}
    First, by \cite[Lemma 5.2]{Keckic2023JMAA}, we see that $(1,0)\strongperp_{BJ}\,(0,A)$ for every $A\in M_2(\mathbb{C})$.

    Next, if $A$ is non-invertible, there exists $A'\in M_2(\mathbb{C})$ such that $A\strongperp_{BJ}\,A'$ (see \cite[Proposition 2.4]{Keckic2023JMAA}), so 
    $$(1,0)\strongperp_{BJ}\,(0,A')\strongperp_{BJ}\,(1,A),$$
    again by \cite[Lemma 5.2]{Keckic2023JMAA}, thus we conclude $e((1,0))=2$.
\end{proof}

Next, we will see that these are not the only vertices in the center of connected component of $\C\oplus M_2(\C)$. 

\begin{lemma}
    Let $A\in M_2(\C)$ be an isometry. Then $e((0,A))=2$.
\end{lemma}
\begin{proof}
    Again, $(0,A)\strongperp_{BJ}(1,0)$ by \cite[Lemma 5.2]{Keckic2023JMAA}. Also, for arbitrary non-zero $B\in M_2(\C)$, by the same argument, it holds 
    $$(0,A)\strongperp_{BJ}(1,0)\strongperp_{BJ}(0,B).$$

    Finally, for arbitrary rank-one $B\in M_2(\C)$ we will have two cases: when $\|B\|\leq1$ and when $\|B\|>1$. In case where $\|B\|\leq1$, we automatically have 
    $$(0,A)\strongperp_{BJ}(1,B).$$
    Namely, by Theorem \ref{glavna}, $A\strongp_{BJ}B$ because $AM_0(A)\cap(\Image B)^{\perp}\neq\{0\}$ from the fact that $AM_0(A)=\C^2$ and $\dim(\Image B)^{\perp}=1$. So by Lemma \ref{C+M_2} we conclude $(0,A)\strongp_{BJ}(1,B)$. The other direction is also true by the same lemma, because $\|(1,B)\|=1$.
    
    In the latter case, when $\|B\|>1$, there is a rank-one $B'\in M_2(\C)$ such that $B\strongperp_{BJ}B'$ and $\|B'\|=1$. Then, according to the previous consideration, the following chain of orthogonalities holds 
    $$(0,A)\strongperp_{BJ}(1,B')\strongperp_{BJ}(1,B),$$
    where first orthogonality is proven in the first case, and the second is by the fact that $B\strongperp_{BJ}B'$ and $\|(1,B')\|=\|B'\|$ and $\|(1,B)\|=\|B\|$. Therefore, $e((0,A))\leq2$ and together with Lemma \ref{ecc>1 for all} we get $e((0,A))=2$.
\end{proof}

We next show that an eccentricity equal to 3 is also possible in $\C\oplus M_2(\C)$.

\begin{lemma}
    Let $A\in M_2(\C)$ be a rank-one operator or a rank-two operator such that $\dim(AM_0(A))=1$. Then $e((0,A))=3$.
\end{lemma}
\begin{proof}
    We can assume $\|A\|=1$ and that 
    $$A=x_1y_1^*+\sigma_2x_2y_2^*,\quad 1>\sigma_2\geq0,$$
    is its SVD so $AM_0(A)=\Span\{x_1\}$. We will first show that every other non-isolated vertex is at distance less than or equal to 3. 

    First, $(0,A)\strongperp_{BJ}(1,0)$ by \cite[Lemma 5.2]{Keckic2023JMAA}.

    Also, $(0,A)\strongperp_{BJ}(1,0)\strongperp_{BJ}(0,B)$ for an arbitrary non-zero $B$ by the previous step.

    Finally, if we consider the vertex $(1,B)$, where $B$ is non-invertible, then there exists non-zero, non-invertible $B'$ such that $B'\strongperp_{BJ}B$. Then, by the same argument and \cite[Lemma 5.2]{Keckic2023JMAA}
    $$(0,A)\strongperp_{BJ}(1,0)\strongperp_{BJ}(0,B')\strongperp_{BJ}(1,B),$$
    so $e((0,A))\leq3$.

    Now, define $$B=\left(\frac{x_1+x_2}{\sqrt{2}}\right)\left(\frac{x_1+x_2}{\sqrt{2}}\right)^*$$ and denote $b=(2,B)$. We will prove that the distance between $a=(0,A)$ and $b$ is equal to 3. First, $a\not\!\strongp_{BJ}b$ by Lemma \ref{C+M_2}, because $A$ attains its norm only on the scalar multiples of $y_1$ and $B^*Ay_1=\frac{x_1+x_2}{2}\neq0$, so $A\not\!\strongp_{BJ}B$. If there exists non-zero $(\lambda,C)\in\C\oplus M_2(\C)$ such that 
    $$a\strongperp_{BJ}(\lambda,C)\strongperp_{BJ}b,$$
    then $\lambda=0$ by Lemma \ref{C+M_2}, because $(2,B)\strongp_{BJ}(\lambda,C)$ and $2>\|B\|$. But then $A\strongperp_{BJ}C$. If $\sigma_2>0$, then $A$ is invertible, and hence no non-zero matrix $C$ can be mutual strong Birkhoff-James orthogonal to $A$. If $\sigma_2=0$, from $A\strongperp_{BJ}C$ we directly have $C=x_2z_2^*$ for some non-zero vector $z_2\in\C^2$ (see \cite[Lemma 2.2]{ArambasicBJMA2020}). But we also have that $C\strongp_{BJ}B$ which is not true by $\langle x_2,\frac{x_1+x_2}{\sqrt{2}}\rangle=\frac{1}{\sqrt{2}}\neq0$. So we have proven that $e((0,A))=3$.
\end{proof}

There is another family of vertices of eccentricity 3.

\begin{lemma}
    Let $A\in M_2(\C)$ be a rank-one operator such that $\|A\|\leq1$. Then $e((1,A))=3$.
\end{lemma}
\begin{proof}
    We will first show that $e((1,A))\leq3$. Take any $(0,B)\in\C\oplus M_2(\C)$. We have that 
    $$(1,A)\strongperp_{BJ}(0,A')\strongperp_{BJ}(1,0)\strongperp_{BJ}(0,B),$$
    by \cite[Lemma 5.2]{Keckic2023JMAA}, where $A'$ is such that $A'\strongperp_{BJ}A$. 
    
    Then, take any $(1,B)$, where $B$ is non-invertible. If $\|B\|\leq1$, then 
    $$(1,A)\strongperp_{BJ}(0,I)\strongperp_{BJ}(1,B).$$
    In the other case when $\|B\|>1$, there exists $B'\in M_2(\C)$ such that $\|B'\|=1$ and $B'\strongperp_{BJ}B$. By Lemma \ref{prelazak na istu normu}, we have $(1,B')\strongperp_{BJ}(1,B)$. Also, $(0,I)\strongperp_{BJ}(1,C)$ for every non-invertible $C\in M_2(\C)$ satisfying $\|C\|\leq1$ by Lemma \ref{C+M_2}, so we have
    $$(1,A)\strongperp_{BJ}(0,I)\strongperp_{BJ}(1,B')\strongperp_{BJ}(1,B).$$
    Next, let $A=\sigma x_1y_1^*, \sigma\leq1$ be its SVD and let $\{x_1,x_2\}$ be an orthonormal basis for $\C^2$. Define $$B=2\left(\frac{x_1+x_2}{\sqrt{2}}\right)\left(\frac{x_1+x_2}{\sqrt{2}}\right)^*.$$
    We will prove that the distance between $(1,A)$ and $(1,B)$ is equal to 3. 
    
    First, $B\not\!\strongp_{BJ} A$, and $\|B\|=2>1$ so $(1,B)\not\!\strongperp_{BJ}(1,A)$ by Lemma \ref{C+M_2}. Furthermore, if there exists non-zero $(\lambda,C)$ such that 
    $$(1,A)\strongperp_{BJ}(\lambda,C)\strongperp_{BJ}(1,B),$$
    then from $(\lambda,C)\strongp_{BJ}(1,B)$, using Lemma \ref{C+M_2}, we get $\|C\|\geq|\lambda|$ and $C\strongp_{BJ}B$ and from $(1,B)\strongp_{BJ}(\lambda,C)$ we get $B\strongp_{BJ}C$, so after rescaling, we may assume
    $$C=\left(\frac{x_1-x_2}{\sqrt{2}}\right)z^*,$$
    for some unit vector $z\in\C^2$ and $|\lambda|\leq1$ (we know that $C$ is non-invertible by $C\strongperp_{BJ}B$ and it is non-zero from $\|C\|\geq|\lambda|$). Then from $(\lambda,C)\strongp_{BJ}(1,A)$, by Lemma \ref{C+M_2}, we conclude $C\strongp_{BJ}A$ which is not true because $\frac{x_1-x_2}{\sqrt{2}}\not\perp x_1$. Together with the previously established upper bound, this gives $e((1,A))=3$. 
\end{proof}

Finally, we end this subsection with a description of vertices with eccentricity equal to 4 (which is maximal possible, because $\diam (\C\oplus M_2(\C))=4$; see \cite[Example 3.2, Lemma 3.5]{Stefanovic}).

\begin{lemma}
    Let $A\in M_2(\C)$ be a rank-one operator such that $\|A\|>1$. Then $e((1,A))=4$.
\end{lemma}
\begin{proof}
    Let $$A=\sigma x_1y_1^*,\quad \sigma>1$$ be its SVD and let $\{x_1,x_2\}$ be an orthonormal basis for $\C^2$. Denote $a=(1,A)$ and define $b=(1,B)$ where   $$B=2\left(\frac{x_1+x_2}{\sqrt{2}}\right)\left(\frac{x_1+x_2}{\sqrt{2}}\right)^*,$$
    so both attain norm on the second component. We will prove that $d(a,b)=4$.
    
    First, $a\not\!\strongp_{BJ} b$ by Lemma \ref{C+M_2} because $A\not\!\strongp_{BJ}B$.

    Next, if $a\strongperp_{BJ}c=(\lambda,C)$, then $A\strongp_{BJ} C$ according to Lemma \ref{C+M_2}. So, by Theorem \ref{glavna}, $AM_0(A)\cap(\Image C)^{\perp}\neq\{0\}$, so $$x_1\in (\Image C)^{\perp}.$$
    Case $C=0$ is not possible again by Lemma \ref{C+M_2}, so $C$ is a rank-one operator and we can assume $\|C\|=1$. Then, its SVD has the form
    $$C=x_2t_2^*$$
    for some unit vector $t_2\in\C^2$. But then $C$ attains its norm only on scalar multiples of $t_2$ and 
    $$B^*Ct_2=x_1+x_2\neq0,$$
    so $C\not\!\strongp_{BJ}B$. But then $c\not\!\strongp_{BJ}b$ by Lemma \ref{C+M_2}. So, there is no non-zero $c\in\C\oplus M_2(\C)$ such that $a\strongperp_{BJ}c\strongperp_{BJ}b$.

    Moreover, from the previous consideration, if it holds 
    $$a\strongperp_{BJ}c\strongperp_{BJ}d\strongperp_{BJ}b,$$
    where $c=(\lambda,C)$ and $d=(\mu,D)$, then $C=x_2t_2^*$. Using the same arguments and the fact that $\|B\|=2>1$, we can conclude that (assuming that also $\|D\|=1$)
    $$D=\left(\frac{x_1-x_2}{\sqrt{2}}\right)z_2^*,$$
    for some unit vector $z_2\in\C^2$. From $\left\langle\frac{x_1-x_2}{\sqrt{2}},x_2\right\rangle=-\frac{1}{\sqrt{2}}\neq0$ we conclude that $C\not\!\strongp_{BJ}D$ and also $D\not\!\strongp_{BJ}C$. From $(\lambda,C)\strongp_{BJ}(\mu,D)$ and $C\not\!\strongp_{BJ}D$, by Lemma \ref{C+M_2}, we get $\mu=0$. Similarly, from $(\mu,D)\strongp_{BJ}(\lambda,C)$ and $D\not\!\strongp_{BJ}C$, by Lemma \ref{C+M_2}, we get $\lambda=0$. So, in that case, $(0,C)\strongperp_{BJ}(0,D)$ would have to hold, which is not true. Thus, we have shown that $e((1,A))=4$ for all non-invertible $A\in M_2(\C),\|A\|>1$.
\end{proof}

After all these lemmas, we proved the following theorem.

\begin{theorem}\label{classification C+M_2}
    Let $a=(\lambda,A)\in\C\oplus M_2(\C)$ be a non-zero, non-invertible vertex. Then    $$e(a)=
\begin{cases}
2,
&
A=0\ \text{or}\ (\lambda=0\ \text{and}\ A\ \text{is a scalar multiple of a unitary} ),\\
4,
& 
\lambda\neq0, \rank A=1\ \text{and}\ \|A\|>|\lambda|, \\
3,
&
\text{in all other cases}.
\end{cases}
$$
\end{theorem}

\subsection{General case}

In this part, we consider all remaining finite-dimensional $C^\ast$-algebras. Thus, throughout this section, $\mathfrak{A}\cong \bigoplus_{i=1}^{k} M_{n_i}(\C), k\geq2$ is a $C^\ast$-algebra that is not isomorphic to $\C\oplus\C$ or $\C\oplus M_2(\C)$, which we have already described. From inclusion
$$\bigoplus_{i=1}^{k} M_{n_i}(\C)\subseteq B\left(\bigoplus_{i=1}^{k}\C^{n_i}\right),$$
and we can see $a$ as an operator on $\bigoplus_{i=1}^{k}\C^{n_i}$. Then from \cite[Lemma 2.1]{KuzmaSushilDebrecen} we know that 
$$M_0(a)=\bigoplus_{i=1}^{k} M_0^*(a_i)$$
where
$$
M_0^*(a_i)=
\begin{cases}
M_0(a_i), & \text{if } \|a_i\|=\|a\|,\\
\{0\}, & \text{if } \|a_i\|<\|a\|.
\end{cases},
$$ 
where $\{0\}$ denotes the zero subspace of $\C^{n_i}$. Also, we will denote 
$$I(a)=\{i\mid \|a_i\|=\|a\|\}.$$

We are now in a position to prove the corresponding variant of the strong Stampfli-Magajna-Bhatia-\v{S}emrl theorem for finite-dimensional $C^\ast$-algebras.

\begin{theorem}\label{SMBS for finite}
    Let $a=(a_1,\dots,a_k),b=(b_1,\dots,b_k)\in\mathfrak{A}\setminus\{0\}$. Then the following conditions are equivalent:
    \begin{itemize}
        \item[i)] $a\strongp_{BJ}b$;
        \item[ii)] there exists $i\in I(a)$ such that $a_iM_0(a_i)\cap(\Image b_i)^{\perp}\neq\{0\}$;
        \item[iii)] $aM_0(a)\cap(\Image b)^{\perp}\neq\{0\}$. 
    \end{itemize}
\end{theorem}
\begin{proof}
    We will first prove ii)$\Rightarrow$i). Suppose that ii) is true. Then for every $c_i\in M_{n_i}(\C)$ 
    $$\|a_i+b_ic_i\|\geq\|a_i\|.$$
    Take an arbitrary $c=(c_1,\dots,c_k)\in\mathfrak{A}$. Then, by the preceding argument and the fact $i\in I(a)$ (so $\|a_i\|=\|a\|$), we have 
    $$\|a+bc\|=\max_{1\leq j\leq k}\|a_j+b_jc_j\|\geq\|a_i+b_ic_i\|\geq\|a_i\|=\|a\|,$$
    so $a\strongp_{BJ}b$.

    Next, assume that ii) is not true. So, for every $i\in I(a)$, we have $a_iM_0(a_i)\cap(\Image b_i)^{\perp}=\{0\}$. By Theorem \ref{glavna}, for every $i\in I(a)$, we can find non-zero $d_i\in M_{n_i}(\C)$ such that $\|a_i+b_id_i\|<\|a_i\|=\|a\|$. Define $c=(c_1,\dots,c_k)$ such that $c_i=d_i$ for every $i\in I(a)$ and $c_i=0$ if $i\notin I(a)$. Then for every $j$ it holds $\|a_j+b_jc_j\|<\|a\|$ so  
    $$\|a+bc\|=\max_{1\leq j\leq k}\|a_j+b_jc_j\|<\|a\|,$$
    so $a\not\!\strongp_{BJ}b$ and we have proven that i) is also not true.

    Finally, ii) and iii) are equivalent because by the consideration preceding the theorem, we have $aM_0(a)=\bigoplus_{i\in I(a)}a_iM_0(a_i)$ and $(\Image b)^{\perp}=\bigoplus_{i=1}^{k}(\Image b_i)^{\perp}$, so 
    $$aM_0(a)\cap(\Image b)^{\perp}=\bigoplus_{i\in I(a)}(a_iM_0(a_i)\cap (\Image b_i)^\perp),$$
    from which the conclusion follows directly.
\end{proof}

We will begin by considering non-isolated vertices that attain the norm in one-dimensional space.

\begin{lemma}\label{jednodimenzionalan u finite}
    Let $a=(a_1,\dots,a_k)\in\mathfrak{A}$ satisfy $\dim (\Image a)^{\perp}\geq1$ and $\dim (aM_0(a))=1$. If $a$ attains the norm in the summand $M_{n_i}(\C),\ n_i\geq2$, then $e(a)=3$.
\end{lemma}
\begin{proof}
    From the assumption, we can assume that $1=\|a_1\|>\|a_i\|$ for all $i\geq2$ and $\dim(a_1M_0(a_1))=1$, so $a_1$ has SVD
    $$a_1=x_1y_1^*+\sum_{i=2}^{n_1}\sigma_ix_iy_i^*,\quad 1>\sigma_2\geq\dots\geq\sigma_{n_1}\geq0.$$
    Define $\displaystyle b=\left(b_1,\frac{1}{2}I_{n_2},\dots,\frac{1}{2}I_{n_k}\right)$, where 
    $$b_1=\frac{x_1+x_{n_1}}{\sqrt{2}}\left(\frac{x_1+x_{n_1}}{\sqrt{2}}\right)^*+\frac{1}{2}\sum_{i=2}^{n_1-1}x_ix_i^*.$$
    If non-zero $c:=(c_1,\dots,c_k)$ holds $c\strongperp_{BJ}b$, then we have $c_1\strongperp_{BJ}b_1$. Indeed, $(\Image b_i)^{\perp}=\{0\}$ for all $i\geq2$, so $c_iM_0(c_i)\cap(\Image b_i)^{\perp}=\{0\}$ for every $i\geq2$. By Theorem \ref{SMBS for finite}, we have $c_1\strongp_{BJ}b_1$  and $\|c_1\|=\|c\|>0$. Also, $\|b_1\|=1>\|b_i\|$ for all $i\geq2$, so $I(b)=\{1\}$. Again, by Theorem \ref{SMBS for finite}, $b_1\strongp_{BJ}c_1$. From this consideration, we know $a\not\!\strongperp_{BJ}b$ ($a_1M_0(a_1)=\Span\{x_1\}$ and $(\Image b_1)^{\perp}=\Span\{\frac{x_1-x_{n_1}}{\sqrt{2}}\}$ so their intersection is $\{0\}$). Finally, if there exists $c=(c_1,\dots,c_k)$ such that $a\strongperp_{BJ}c\strongperp_{BJ}b$, then $c_1\strongperp_{BJ}b_1$ by the preceding argument and also $a_1\strongp_{BJ}c_1$ by the assumption on $a$ ($I(a)=\{1\}$). Thus, $a_1\strongp_{BJ}c_1\strongp_{BJ}b_1,$ which contradicts Remark \ref{remark:one-sided-AM0}. Hence, $e(a)\geq3$ and since the diameter of the non-trivial connected component is 3, it follows that $e(a)=3$.
\end{proof}

Now we can move on to studying the case when $a$ attains the norm only on one summand, which is $\C$.

\begin{lemma}\label{centralni skalarni cvor}
    If $\mathfrak{A}\cong\C\oplus\mathfrak{A}'$ is a finite-dimensional $C^\ast$-algebra, then the eccentricity of the vertex $a=(1,0,\dots,0)$ is 2. 
\end{lemma}
\begin{proof}
    Choose an arbitrary non-zero non-invertible vertex $b=(b_1,b_2,\dots,b_k)$. If $b_1=0$, then $a\strongperp_{BJ}b$ by \cite[Lemma 5.2.]{Keckic2023JMAA}. If $b_1\neq0$ (so $a\not\!\strongperp_{BJ}b$), then there exists some $b_{\ell}\in M_{n_l}$ which is non-invertible element. So, we can find $b_{\ell}'\neq0$ such that $b_{\ell}'\strongperp_{BJ}b_{\ell}$ (in the case where $n_l=1$, we know that $b_l=0$, so we can choose any non-zero $b_l'$). Define $c$ so that $c_{\ell}=b_{\ell}'$ is the only non-zero element of $c$. Then, directly by \cite[Lemma 5.2]{Keckic2023JMAA}, we conclude that $a\strongperp_{BJ}c\strongperp_{BJ}b$, so $e(a)\leq2$. Lemma \ref{ecc>1 for all} yields $e(a)\geq2$, therefore $e(a)=2$.
\end{proof}

\begin{lemma}\label{jednodimenzionalan skalarni}
    Let $\mathfrak{A}\cong\C\oplus\mathfrak{A}'$ be a finite-dimensional $C^\ast$-algebra. If $a\in\mathfrak{A}$ satisfies $\dim(\Image a)^{\perp}\geq1$ and $\dim (aM_0(a))=1$, where $a$ attains the norm in the summand $\C$ and $[a]\neq[(1,0,\dots,0)]$, then the eccentricity of $a$ is equal to 3.
\end{lemma}
\begin{proof}
    After rescaling, we can assume that $a=(1,a'), a'\in\mathfrak{A}'$ is such that $0<\|a'\|<1$. Define $b=(0,b')$ so that $b'M_0(b')=\Image a'$ and $b'=(b_1',\dots,b_{k-1}')$ is invertible. Such $b'$ exists. Namely, if $a'=(a_1',\dots,a_{k-1}')$, let $E_i=\Image a_i', 1\leq i\leq k-1$. Then we can define 
    \begin{align*}
    b_i'=\left\{
    \begin{array}{cc}
    P_{E_i}+\frac{1}{2}P_{E_i^{\perp}}, & \text{if }E_i\neq\{0\},\\
    \frac{1}{2}I, & \text{if } E_i=\{0\}.
    \end{array}\right.
    \end{align*}
    Then $b\not\!\strongp_{BJ}a$ because $b'M_0(b')\cap(\Image a')^{\perp}=\{0\}$. Also, if 
    $$(1,a')\strongperp_{BJ}(\lambda,c')\strongperp_{BJ}(0,b'),$$
    then $\lambda$ must be zero from $(1,a')\strongperp_{BJ}(\lambda,c')$ and $1>\|a'\|$. But from $(0,c')\strongperp_{BJ}(0,b')$ we get $c'\strongperp_{BJ}b'$ which is not possible by the invertibility of $b'$. Thus $e(a)\geq3$ and since the diameter of the non-trivial component is 3, it follows that $e(a)=3$.
\end{proof}

Now we will move on to the case when $\dim(\Image a)^{\perp}=1$. First, we will consider the situation when the only non-invertible component belongs to block $M_n(\C), n\geq2$.

\begin{lemma}\label{kodimenzija jedan nekomutativni}
    If $a\in\mathfrak{A}$ is such that $\dim(\Image a)^{\perp}=1$, where the unique non-invertible block of $a$ belongs to a noncommutative summand, then $e(a)=3$.
\end{lemma}
\begin{proof}
    We can assume that $a=(a_1,a_2,\dots,a_k)\in\mathfrak{A}$ is such that $\|a_1\|=1,\rank a_1=n_1-1$ with SVD 
    $$a_1=\sum_{i=1}^{n_1-1}\sigma_ix_iy_i^*,\ 1=\sigma_1\geq\dots\geq\sigma_{n_1-1}>0$$
    and $a_i$ is invertible for $i\geq2$. Let $\{x_1,\dots,x_{n_1-1},x_{n_1}\}$ be an orthonormal basis for $\C^{n_1}$ and define $\displaystyle b=\left(b_1,\frac{1}{2}I,\dots,\frac{1}{2}I\right)$, where 
    $$b_1=\frac{x_1+x_{n_1}}{\sqrt{2}}\left(\frac{x_1+x_{n_1}}{\sqrt{2}}\right)^*.$$  
    Thus $b\not\!\strongp_{BJ}a$ because $I(b)=\{1\}$ and 
    $$b_1M_0(b_1)\cap(\Image a_1)^{\perp}=\{0\}.$$ 
    Also if for non-zero $c\in\mathfrak{A}$ it holds $a\strongperp_{BJ}c\strongperp_{BJ}b$, then $c_1\strongperp_{BJ}b_1$ and $\|c_1\|=\|c\|>0$ as in the proof of Lemma \ref{jednodimenzionalan u finite}. Finally, from $c\strongp_{BJ}a$ and the hypothesis of the lemma, we get 
    $$\max\{\|c_1+a_1x_1\|,\dots,\|c_k+a_kx_k\|\}\geq\max\{\|c_1\|,\dots,\|c_k\|\}\geq\|c_1\|, \forall (x_1,\dots,x_k)\in\mathfrak{A},$$
    so if we put $x_i=-a_i^{-1}c_i$ for $i\geq2$, we get $c_1\strongp_{BJ}a_1$. So we have a chain $b_1\strongp_{BJ}c_1\strongp_{BJ}a_1$, which contradicts Remark \ref{remark:one-sided-kernel}. Thus $e(a)\geq3$ and since the diameter of the non-trivial component is 3, it follows that $e(a)=3$.
\end{proof}

In the next lemma, we will solve the remaining cases when $\dim (\Image a)^{\perp}=1$.

\begin{lemma}\label{kodimenzija jedan skalarni}
    Let $a\in\mathfrak{A}\cong\C\oplus\mathfrak{A}'$ be such that $\dim(\Image a)^{\perp}=1$, where the unique non-invertible block of $a$ belongs to $\C$. If $aM_0(a)=\{0\}\oplus\bigoplus_{i=2}^{k}\C^{n_i}$, then $e(a)=2$; otherwise $e(a)=3$.
\end{lemma}
\begin{proof}
    We can assume 
    $$a=(0,a_2,\dots,a_k)=(0,a'),$$
    such that $\|a\|=1$, where all $a_i,i\geq2$ are invertible. If $aM_0(a)\neq\{0\}\oplus\bigoplus_{i=2}^{k}\C^{n_i}$, then we can find some summand (let us say $a_2$) such that $a_2$ is not a unitary, or $\|a_2\|<\|a\|$. If $\|a_2\|<\|a\|$, choose an arbitrary unit vector  $v\in\C^{n_2}$; if $\|a_2\|=\|a\|$, but $a_2$ is not a unitary, we can find a unit vector $v\not\in a_2M_0(a_2)$. In both cases, we can define $b$ as
    $$b=\left(1,\frac{1}{2}P_{\Span\{v\}^{\perp}},\frac{1}{2}I,\dots,\frac{1}{2}I\right),$$
    which is the non-isolated vertex such that $a\not\!\strongp_{BJ}b$ by Theorem \ref{SMBS for finite} because $aM_0(a)\cap(\Image b)^{\perp}=\{0\}$. 

    Also, if there exists non-zero $c=(c_1,c')$ such that $a\strongperp_{BJ}c\strongperp_{BJ}b$, then from $b\strongp_{BJ}c$ and the fact that $\|b_1\|>\|b_i\|,i\geq2$, we get $c_1=0$. But then $c'\strongperp_{BJ}a'$ which is not possible because $a'$ is invertible. Thus, we have proven that in this case $e(a)\geq3$ and since the diameter of the non-trivial connected component is 3, it follows that $e(a)=3$.

    In the remaining case when $aM_0(a)=\{0\}\oplus\bigoplus_{i=2}^{k}\C^{n_i}$, we know that
    $$a=(0,u_2,\dots,u_k),$$
    where all the matrices $u_i$ are unitary. Take an arbitrary non-zero non-invertible $b=(b_1,b_2,\dots,b_k)$. If $b_1=0$, then $d(a,b)\leq2$ because
    $$a\strongperp_{BJ}(1,0,\dots,0)\strongperp_{BJ}b,$$
    by \cite[Lemma 5.2]{Keckic2023JMAA}. In what follows, we may therefore assume that $b_1\neq0$. Then there exists an index $p$, with $2\leq p\leq k$, such that $b_p$ is non-invertible. If $1\in I(b)$ then $b\strongp_{BJ}a$ by Theorem \ref{SMBS for finite}. Also, $a\strongp_{BJ}b$ by the same theorem because $a_pM_0(a_p)\cap(\Image b_p)^{\perp}\neq\{0\}$, so $a\strongperp_{BJ}b$. So, in the remaining case, there is $q$ such that $2\leq q\leq k$ and $q \in I(b)$. Define $c=(c_1,\dots,c_k)$ such that only non-zero coordinates are $c_p=b_p'$ and $c_1=1$ where $b_p'$ is  element such that $b_p'\strongperp_{BJ}b_p$ and $\|b_p'\|=1$. Such $b_p'$ does exist - if $n_p>1$ it can be chosen to be non-invertible, and if $n_p=1$ (then $b_p=0$), we can choose $b_p'=1$. Then $1\in I(c)$ and $a_1=0$ so by Theorem \ref{SMBS for finite} we have $c\strongp_{BJ}a$. From the same theorem and facts that $p\in I(c)$ and $c_p\strongp_{BJ}b_p$ we get $c\strongp_{BJ}b$. For the other direction, we first assume $p\neq q$. Then $q\in I(a)$ and $c_q=0$ so $a\strongp_{BJ}c$; also $q\in I(b)$ so $b\strongp_{BJ}c$. In the case when $p=q$, from $p=q\in I(b)$, we have $b_p\neq0$, so from $b_p\strongp_{BJ}b_p'$ it follows that $b_p'$ is non-invertible so $(\Image b_p')^{\perp}\neq\{0\}$. Then $p\in I(a)\cap I(b)$ and $$a_pM_0(a_p)\cap(\Image b_p')^{\perp}\neq\{0\}\ \text{and}\ b_pM_0(b_p)\cap(\Image b_p')^{\perp}\neq\{0\},$$ so $a\strongp_{BJ}c$ and $b\strongp_{BJ}c$, respectively, by Theorem \ref{SMBS for finite}. Thus, we have proven that $e(a)\leq2$ and by Lemma \ref{ecc>1 for all} we get $e(a)=2$. 
\end{proof}

\begin{lemma}\label{oba prostora bar dva}
    If $a\in\mathfrak{A}$ is such that $\dim(\Image a)^{\perp}\geq2$ and $\dim(aM_0(a))\geq2$, then $e(a)=2$.
\end{lemma}
\begin{proof}
    The core of the proof is contained in Lemma \ref{oba bar dva}. Take any non-zero non-invertible $b=(b_1,\dots,b_k)$. Let $p\in I(b)$ and let $u\in bM_0(b)$ be a unit vector such that $u_i=0$ for every $i\neq p$. Also, let $q$ be an index such that $b_q$ is non-invertible and let $t\in(\Image b)^{\perp}$ be a unit vector such that $t_i=0$ for every $i\neq q$. Furthermore, we can choose a unit vector 
    $$v\in aM_0(a)\cap t^{\perp}$$
    by a dimension argument ($\dim(aM_0(a))\geq2$). Also, we can choose a unit vector
    $$w\in (\Image a)^{\perp}\cap u^{\perp}$$
    by the same argument ($\dim(\Image a)^{\perp}\geq2$). Finally, define $c:=(c_1,\dots,c_k)\in\mathfrak{A}$ such that $$c_i:=P_{\Span\{t_i,w_i\}}$$ is an orthogonal projection onto a linear span of the $i$-th coordinates of vectors $t$ and $w$. 

    First, see that $v_i\perp t_i$ for every $i$. Namely, $t_i=0$ for every $i\neq q$ and from $v\in t^{\perp}$, we get $v_q\perp t_q$. For the same reason, we get $w_i\perp u_i$ for every $i$. Next, from $v\in aM_0(a)$ we get $v=ax$ (so $v_i=a_ix_i$ for all $1\leq i\leq k$) for some $x=(x_1,\dots,x_k)\in M_0(a)$. If $\|a_i\|=\|a\|$, then $x_i\in M_0(a_i)$ and so $v_i\in a_iM_0(a_i)\subseteq\Image a_i$; if $\|a_i\|<\|a\|$, then $x_i=0$ and so $v_i=0$. On the other hand, $w\in (\Image a)^{\perp}=\bigoplus_{i=1}^{k}(\Image a_i)^{\perp}$, so $w_i\in(\Image a_i)^{\perp}$, from which we can conclude $v_i\perp w_i$ for every $1\leq i\leq k$. Finally, $u_i\perp t_i$ for every $i$. Indeed, the only non-zero coordinates are $u$ and $t$ are on the positions $p$ and $q$, respectively. If $p\neq q$, this is immediate. If $p=q$, then $u_p\in\Image b_p$ and $t_p\in(\Image b_p)^{\perp}$ so they are orthogonal. 

    From this consideration, we have 
    $$c_iv_i=c_iu_i=0\quad,\quad c_iw_i=w_i \quad\text{and}\quad c_it_i=t_i.$$ Now, we can apply Theorem \ref{SMBS for finite} four times: $v\in aM_0(a)\cap (\Image c)^{\perp}$ so $a\strongp_{BJ}c$; $w\in cM_0(c)\cap (\Image a)^{\perp}$ so $c\strongp_{BJ}a$; $t\in cM_0(c)\cap (\Image b)^{\perp}$ so $c\strongp_{BJ}b$; $u\in bM_0(b)\cap (\Image c)^{\perp}$ so $b\strongp_{BJ}c$ and we finally get $a\strongperp_{BJ}c\strongperp_{BJ}b$, thus $e(a)\leq2$. Lemma \ref{ecc>1 for all} yields $e(a)\geq2$, therefore $e(a)=2$.
\end{proof}

So we finally have full characterization of eccentricity in finite-dimensional $C^\ast$-algebra.

\begin{theorem}\label{the main}
    Let $\mathfrak{A}\cong\bigoplus_{i=1}^{k}M_{n_i}(\C)$ be a finite-dimensional $C^\ast$-algebra different from $\C,\C\oplus\C$ and $M_2(\C)$. Then all non-zero non-invertible elements belong to the same connected component. Moreover, the eccentricity of a vertex $a=(a_1,\dots,a_k)$ in this component is given by
    $$e(a)=
\begin{cases}
2,
&
\begin{array}{l}
\dim aM_0(a)\geq2\ \text{and}\ \dim (\Image a)^{\perp}\geq2,\ \text{or}\\
a_i=0 \ \text{or}\ a_i=\|a\|u_i \ \text{for all}\ i, \ \text{where each}\ u_i\ \text{is a unitary};
\end{array}
\\[5mm]
4,
&\begin{array}{l}
\mathfrak A\cong M_3(\mathbb C),\quad
\rank a=2,\
\dim aM_0(a)=1,\ \text{or}\\
\mathfrak A\cong\mathbb C\oplus M_2(\mathbb C),
\quad a=(\lambda,A),\ \lambda\neq0,\ \rk A=1,\ \|A\|>|\lambda|;
\end{array}
\\[5mm]
3,
&
\ \ \text{in all other cases}.
\end{cases}
$$
\end{theorem}
\begin{proof}
    The case $\mathfrak{A}\cong M_n(\C)$ follows from Theorems \ref{M_3} and \ref{M_n,n>3}. The exceptional algebra $\C\oplus M_2(\C)$ is covered by Theorem \ref{classification C+M_2}. We may therefore assume that $\mathfrak{A}$ has at least two summands and that $\mathfrak{A}\not\cong\C\oplus M_2(\C).$

    Let $a\in\mathfrak{A}$ be a non-zero non-invertible element. Since $a$ is non-zero and non-invertible, we have $\dim aM_0(a),\dim(\Image a)^{\perp}\geq1$.

    If $\dim aM_0(a)\geq2$ and $\dim(\Image a)^{\perp}\geq2$, then Lemma \ref{oba prostora bar dva} gives $e(a)=2$.

    Suppose next that $\dim aM_0(a)=1$. In this case, $a$ attains its norm in exactly one summand. If this is a noncommutative summand, then Lemma \ref{jednodimenzionalan u finite} gives $e(a)=3$. If the norm is attained in a scalar summand, then, after reordering the summands and rescaling, we may write $a=(1,a'),\ \|a'\|<1.$ If $a'=0$, then Lemma \ref{centralni skalarni cvor} gives $e(a)=2$; otherwise Lemma \ref{jednodimenzionalan skalarni} gives $e(a)=3$.

    It remains to consider the case $\dim aM_0(a)\geq2$ and $\dim(\Image a)^\perp=1$. Then precisely one block of $a$ is non-invertible. If this block belongs to a noncommutative summand, Lemma \ref{kodimenzija jedan nekomutativni} gives $e(a)=3$. If it belongs to a scalar summand, then Lemma \ref{kodimenzija jedan skalarni} shows that $e(a)=2$ precisely when, after reordering the summands, $a=(0,\|a\|u_2,\dots,\|a\|u_k),$ where every $u_i$ is unitary; in all other cases, $e(a)=3$.
\end{proof}
\begin{remark}
    The condition appearing in the eccentricity 2 case can also be written more intrinsically. Namely,
$$ a_i=0 \quad\text{or}\quad a_i=\|a\|u_i \quad\text{for all }i,\ \text{where each}\ u_i\ \text{is a unitary} $$
is equivalent to
$$ a^*a=\|a\|^2p $$
for some central projection \(p\in\mathfrak{A}\). Indeed, in the decomposition $\mathfrak{A}\cong\bigoplus_{i=1}^kM_{n_i}(\C),$ every central projection has the form
$$p=(c_1I_{n_1},\dots,c_kI_{n_k}), \qquad c_i\in\{0,1\}.$$
Thus the second condition in Theorem \ref{the main} may be replaced by the shorter and representation-independent condition.
\end{remark}

\begin{remark}
    For completeness, we again (see also \cite[Proposition 2.4]{ArambasicBJMA2020} and \cite[Proposition 4.1]{Keckic2023JMAA}) describe the orthographs of the three small algebras excluded from Theorem \ref{the main}. 

    The graph $\Gamma(\C)$ consists of a single isolated vertex $[1]$ whose eccentricity is 0. 

    In $\C\oplus\C$, invertible elements are isolated vertex. The unique non-trivial connected component consists of the two vertices $[(1,0)]$ and $[(0,1)]$, and they are mutual strong Birkhoff-James orthogonal. Thus, both have eccentricity 1.

    Finally, in $M_2(\C)$, every non-isolated vertex is represented by a rank-one matrix. By Theorem \ref{glavna}, for two rank-one matrices $A=xy^*$ and $B=zw^*$, where $x,y,z,w\in\C^2$, $A\strongp_{BJ}B$ is equivalent to $\langle x,z\rangle=0$ and this is also equivalent to $B\strongp_{BJ}A$. Thus every non-trivial connected component is represented by some non-zero vector $x\in\C^2$ and has the form 
    $$C_x=\{[xy^*] \mid 0\neq y\in\C^2 \}\cup\{[x'z^*] \mid 0\neq z\in\C^2\},$$
    where $x'\in\C^2$ is any unit vector such that $\langle x',x\rangle=0$. Every vertex in $C_x$ has eccentricity 2.
\end{remark}

\subsection{Examples and consequences of the main theorem}

We now record several consequences of Theorem \ref{the main} and give some examples illustrating the classification. The first among them concerns the radius of the orthograph.

\begin{corollary}\label{radius}
    Let $\mathfrak{A}$ be a finite-dimensional $C^\ast$-algebra not isomorphic to any of $\C$, $\C\oplus\C$, and $M_2(\C)$. Then the radius of the unique non-trivial connected component of $\Gamma(\mathfrak{A})$ is
    $$r(\Gamma(\mathfrak{A}))=
        \begin{cases}
            3, & \mathfrak{A}\cong M_3(\C),\\
            2, & \mathfrak{A}\not\cong M_3(\C).
        \end{cases}
    $$
\end{corollary}

\begin{proof}
    If $\mathfrak{A}\cong M_3(\C)$, then Theorem \ref{M_3} shows that every vertex in the unique non-trivial connected component has eccentricity at least $3$ and that vertices of eccentricity $3$ exist. Hence, the radius is $3$.

    Suppose now that $\mathfrak{A}\not\cong M_3(\C)$. If $\mathfrak{A}\cong\bigoplus_{i=1}^kM_{n_i}(\C), \ k\geq2,$ then $$a=(I_{n_1},0,\dots,0)$$ is a non-zero non-invertible element whose every component is either zero or a unitary. Therefore, Theorem \ref{the main} gives $e(a)=2$, and by Lemma \ref{ecc>1 for all} we get $r(\Gamma(\mathfrak{A}))=2$.

    The only remaining case is $\mathfrak{A}\cong M_n(\C)$ with $n\geq4$. Let $p$ be an orthogonal projection of rank $2$. Then
    $$\dim pM_0(p)=2\quad\text{and}\quad\dim(\Image p)^\perp=n-2\geq2,$$
    so Theorem \ref{the main} again gives $e(p)=2$. Since no vertex has eccentricity 1 by Lemma \ref{e>1}, the radius is again $2$.
\end{proof}

    \begin{corollary}
    Among finite-dimensional $C^\ast$-algebras, $M_3(\C)$ and $\C\oplus M_2(\C)$ are the only ones whose non-trivial connected components have diameter $4$. Although their diameters coincide, their radii are different:
    $$r(\Gamma(M_3(\C)))=3\quad\text{and}\quad r(\Gamma(\C\oplus M_2(\C)))=2.$$
\end{corollary}

The following corollary describes the center of the non-trivial connected component of $\Gamma(\mathfrak A)$, where $\mathfrak A$ is a finite-dimensional $C^\ast$-algebra.

\begin{corollary}
    Let $\mathfrak{A}\cong\bigoplus_{i=1}^k M_{n_i}(\C)$ be a finite-dimensional $C^\ast$-algebra not isomorphic to any of $\C$, $\C\oplus\C$, and $M_2(\C)$. Then the center of the unique non-trivial connected component of $\Gamma(\mathfrak{A})$ is described as follows.

    \begin{enumerate}
        \item If $\mathfrak{A}\cong M_3(\C)$, then the center consists
        precisely of the projective classes of non-zero singular
        matrices $a\in M_3(\C)$ such that either
        $$\rank a=1\quad \text{or}\quad \left(\rank a=2 \ \text{and}\ aM_0(a)=\Image a\right).$$

        \item If $\mathfrak{A}\not\cong M_3(\C)$, then the center
        consists precisely of the projective classes of non-zero non-invertible elements $a=(a_1,\dots,a_k)\in\mathfrak{A}$ satisfying at least one of the following two conditions:
        $$\dim aM_0(a)\geq2\quad\text{and}\quad\dim(\Image a)^\perp\geq2,$$
        or, for every $1\leq i\leq k$,
        $$a_i=0\quad\text{or}\quad a_i=\|a\|u_i$$
        for some unitary $u_i\in M_{n_i}(\C)$.
    \end{enumerate}
\end{corollary}

\begin{proof}
    If $\mathfrak{A}\cong M_3(\C)$, then the radius of the non-trivial component is $3$. Therefore, its center consists precisely of the vertices of eccentricity $3$, whose description is given by Theorem \ref{M_3}.

    Suppose that $\mathfrak{A}\not\cong M_3(\C)$. By Corollary \ref{radius}, the radius of the unique non-trivial component is $2$. Consequently, its center consists precisely of the vertices of eccentricity $2$. Their description now follows directly from Theorem \ref{the main}.
\end{proof}

The following corollary shows how, only on the basis of eccentricity, we can distinguish the scalar multiples of unitary matrices in relation to all other invertible elements.

\begin{corollary}
    Let $\mathfrak{B}$ be a finite-dimensional $C^\ast$-algebra with $\dim\mathfrak{B}\geq2$, and let $b\in\mathfrak{B}$ be an invertible element. Then the vertex $(0,b)\in\C\oplus\mathfrak{B}$ has eccentricity $2$ if and only if $b$ is a scalar multiple of a unitary element of $\mathfrak{B}$.
\end{corollary}

\begin{proof}
    Set $a=(0,b)$. Since $b$ is invertible, $\dim(\Image a)^\perp=1.$ Therefore, the first alternative for eccentricity 2 in Theorem \ref{the main} cannot occur. Consequently, $e(a)=2$ if and only if every non-zero component of $b$ is equal to $\|b\|$ times a unitary. This is equivalent to
    $$b=\|b\|u$$
    for some unitary $u\in\mathfrak{B}$.
\end{proof}

In \cite[Example 3.1]{Stefanovic}, it was shown that the diameter of $\C^k,k\geq3$ is equal to $3$. Here we show that this is, in fact, the typical example of a vertex with eccentricity $3$.

\begin{corollary}
     Let $a=(a_1,\dots,a_k)\in\C^k$, where $k\geq3$, be a non-zero non-invertible element satisfying $\|a\|=1$. Then $e(a)=3$ if and only if there exists at least one index $i_0\in\{1,\dots,k\}$ such that $$\left(0<|a_{i_0}|<1\right)\quad \text{and}\quad  \left(|\{i\mid \|a_i\|=\|a\|\}|=1 \text{ or } |\{i\mid a_i=0\}|=1\right).$$
\end{corollary}
\begin{proof}
    We have $$\dim aM_0(a)=|\{i\mid\|a_i\|=\|a\|\}|\quad\text{and}\quad\dim(\Image a)^\perp=|\{i\mid a_i=0\}|.$$
    Hence, Theorem \ref{the main} shows that $e(a)=2$ precisely when there is no index $i$ such that $0<|a_i|<1$ or
    $|\{i\mid \|a_i\|=\|a\|\}|\geq2$ and $|\{i\mid a_i=0\}|\geq2$. Since every non-isolated vertex of $\Gamma(\C^k)$ has eccentricity either $2$ or $3$, the stated characterization follows.
\end{proof}

\section{Orthograph of $K(H),\dim H=\infty$}

Throughout this section, $H$ denotes a separable infinite-dimensional complex Hilbert space, and $K(H)$ denotes the $C^\ast$-algebra of compact operators on $H$. Here for $A\in K(H)$ we use the same notation
$$M_0(A)=\{x\in H \mid \|Ax\|=\|A\|\|x\|\}.$$
If $A\neq0$, compactness ensures that this subspace is non-zero and finite-dimensional.

We first recall the description of the isolated vertices of $\Gamma(K(H))$. By \cite[Example 2.2]{Keckic2023JMAA}, a non-zero operator $A\in K(H)$ represents an isolated vertex if and only if
$$\overline{\Image A}=H \Longleftrightarrow\ker A^*=\{0\}.$$

We shall use the following compact-operator version of the strong Stampfli, Magajna, Bhatia--\v{S}emrl theorem; see \cite[Proposition 2.8]{ArambasicAFA2014} for the corresponding sequential characterization and \cite[Lemma 4.1]{KST} for the formulation below. 

\begin{lemma}\label{Bhatia-Semrl compact}
Let $A,B\in K(H)$, with $A\neq0$. Then
$$A\strongp_{BJ}B\quad\Longleftrightarrow\quad AM_0(A)\cap\ker B^*\neq\{0\}.$$
Equivalently, $A\strongp_{BJ}B$ if and only if there exists a unit vector $x\in M_0(A)$ such that $B^*Ax=0$.
\end{lemma}

In the following theorem, we prove that the unique non-trivial connected component of $\Gamma(K(H))$ has a diameter 3. At the same time, we extend the corresponding result of Lemma \ref{AM_0(A) jednodimenzionalan} by showing that every non-isolated $A\in K(H)$ satisfying $\dim AM_0(A)=1$ has eccentricity $3$.

\begin{theorem}
Let $H$ be an infinite-dimensional separable Hilbert space. Then the unique non-trivial connected component of $\Gamma(K(H))$ has diameter 3. Moreover, every non-isolated vertex $A\in K(H)$ satisfying $\dim AM_0(A)=1$ has eccentricity $3$.
\end{theorem}
\begin{proof}
Let $A,B\in K(H)$ be two non-isolated vertices. Choose unit vectors
$$x\in AM_0(A),\quad y\in BM_0(B),\quad z\in\ker A^*\quad\text{and}\quad t\in\ker B^*.$$
Since $x\in\Image A$ and $z\in(\Image A)^\perp$, we have $x\perp z$. Similarly, $y\perp t$.

Choose two orthonormal vectors
$$u,v\in\Span\{x,y,z,t\}^{\perp},$$
and let $P$ and $Q$ be the orthogonal projections onto 
$$\Span\{z,u\}\quad\text{and}\quad\Span\{t,v\},$$
respectively. Both $P$ and $Q$ have rank two and therefore belong to $K(H)$.

We have $x\in AM_0(A)\cap\ker P\ \text{and}\ z\in PM_0(P)\cap\ker A^*,$ and therefore $A\strongperp_{BJ}P$ by Lemma \ref{Bhatia-Semrl compact}. Similarly, $u\in PM_0(P)\cap\ker Q\ \text{and}\ v\in QM_0(Q)\cap\ker P,$ so $P\strongperp_{BJ}Q$ by the same lemma. Finally, $t\in QM_0(Q)\cap\ker B^*\ \text{and}\ y\in BM_0(B)\cap\ker Q,$ which gives $Q\strongperp_{BJ}B$. Consequently, $A\strongperp_{BJ}P\strongperp_{BJ}Q\strongperp_{BJ}B,$ and hence the diameter is at most $3$.

It remains to show that there are two vertices at distance 3. In fact, we will prove the second statement of the theorem, from which the desired conclusion follows immediately. Let $A\in K(H)$ be a non-isolated vertex such that
$$AM_0(A)=\Span\{x_1\},$$
where $x_1$ is a unit vector. Since $A$ is non-isolated, we have $\ker A^*\neq\{0\}$, so we can choose a unit vector $x_0\in\ker A^*$. Since $x_1\in\Image A$, it follows that $x_1\perp x_0$. Then also $\frac{x_1+x_0}{\sqrt{2}}$ and $\frac{x_1-x_0}{\sqrt{2}}$ are orthonormal vectors. Then extend $\left\{\frac{x_1+x_0}{\sqrt{2}},\frac{x_1-x_0}{\sqrt{2}}\right\}$ to an orthonormal basis $\left\{\frac{x_1+x_0}{\sqrt{2}},\frac{x_1-x_0}{\sqrt{2}},x_2,x_3,\dots\right\}$ of $H$. Define
$$B=\frac{x_1+x_0}{\sqrt{2}}\left(\frac{x_1+x_0}{\sqrt{2}}\right)^*+\sum_{n=2}^{\infty}\frac{1}{n}x_nx_n^*.$$
Then $B$ is a compact selfadjoint operator satisfying
$$BM_0(B)=\Span\left\{\frac{x_1+x_0}{\sqrt{2}}\right\}\quad\text{and}\quad\ker B^*=\Span\left\{\frac{x_1-x_0}{\sqrt{2}}\right\}.$$
In particular, $B$ is also a non-isolated vertex. First, $A\not\!\strongp_{BJ}B$ because
$$AM_0(A)\cap\ker B^*=\Span\{x_1\}\cap\Span\left\{\frac{x_1-x_0}{\sqrt{2}}\right\}=\{0\}.$$

Suppose that there exists a non-zero $C\in K(H)$ such that $A\strongperp_{BJ}C\strongperp_{BJ}B$. From $A\strongp_{BJ}C$ and the fact that $AM_0(A)=\Span\{x_1\}$, we obtain $x_1\in\ker C^*$ by Lemma \ref{Bhatia-Semrl compact}. On the other hand, $C\strongp_{BJ}B$ implies $CM_0(C)\cap\ker B^*\neq\{0\}.$ Since $\ker B^*=\Span\{\frac{x_1-x_0}{\sqrt{2}}\}$, it follows that $\frac{x_1-x_0}{\sqrt{2}}\in CM_0(C)\subseteq\Image C.$ But this is impossible because $x_1\in\ker C^*=(\Image C)^\perp$, but $\langle x_1,\frac{x_1-x_0}{\sqrt{2}}\rangle=\frac{1}{\sqrt{2}}\neq0.$ Consequently, $d(A,B)\geq3$. Since the diameter is at most $3$, we conclude that $e(A)=3$, so the statement holds.
\end{proof}

We now extend Lemma \ref{ImA n-1-dimenzionalan} to compact operators, obtaining another family of vertices with eccentricity $3$.

\begin{lemma}
Let $A\in K(H)$ be a non-zero operator such that $\dim\ker A^*=1$. Then $e(A)=3$.
\end{lemma}
\begin{proof}
Choose unit vectors $x_0\in\ker A^*$ and $x_1\in AM_0(A)$. Since $x_0\perp x_1$, the operator
$$B=\frac{x_0+x_1}{\sqrt{2}}\left(\frac{x_0+x_1}{\sqrt{2}}\right)^*$$
is a rank-one orthogonal projection and so belongs to $K(H)$. The same argument as in the proof of Lemma \ref{ImA n-1-dimenzionalan} shows that $A\not\!\strongperp_{BJ}B$ and that there is no non-zero $C\in K(H)$ such that $A\strongperp_{BJ}C\strongperp_{BJ}B$. Thus, $d(A,B)\geq3$, while the preceding theorem gives $e(A)\leq3$. Consequently, $e(A)=3$.
\end{proof}

Finally, we extend Lemma \ref{oba bar dva} to compact operators and show that all remaining non-isolated vertices have eccentricity $2$.

\begin{lemma}
Let $A\in K(H)$ be a non-zero operator satisfying
$$\dim AM_0(A)\geq2\quad\text{and}\quad\dim\ker A^*\geq2.$$
Then $e(A)=2$.
\end{lemma}
\begin{proof}
The construction in the proof of Lemma \ref{oba bar dva} applies here as well, with $(\Image A)^\perp=\ker A^*$ and $(\Image B)^\perp=\ker B^*$. The operator $C$ satisfying $A\strongperp_{BJ}C\strongperp_{BJ}B$ is an orthogonal projection of rank at most two and is therefore compact. Thus $e(A)\leq2$.

On the other hand, a vertex with these properties cannot have eccentricity 1. Indeed, we can choose a unit vector $x\in AM_0(A)$ and consider rank-one orthogonal projection $P=xx^*$. Thus, $P$ is a non-isolated vertex ($\ker P\neq\{0\}$), $[P]\neq[A]$ and satisfies $P\not\!\strongp_{BJ}A$, since
$$PM_0(P)\cap\ker A^*=\Span\{x\}\cap\ker A^*=\{0\}.$$
Therefore, $e(A)=2$.
\end{proof}

We can now state a theorem that describes the eccentricity of an arbitrary non-isolated vertex in $\Gamma(K(H))$. It also determines the radius and describes the center of the unique non-trivial connected component.

\begin{theorem}
Let $H$ be an infinite-dimensional separable complex Hilbert space. The unique non-trivial connected component of $\Gamma(K(H))$ has a diameter $3$ and a radius $2$. For every non-isolated vertex $A\in K(H)$, its eccentricity within this component is given by
$$e(A)=\begin{cases}
3, & \dim AM_0(A)=1 \text{ or } \dim\ker A^*=1,\\
2, & \dim AM_0(A)\geq2 \text{ and } \dim\ker A^*\geq2.
\end{cases}
$$
Consequently, the center of this component consists precisely of the projective classes of non-zero operators $A\in K(H)$ satisfying
$$\dim AM_0(A)\geq2\qquad\text{and}\qquad\dim\ker A^*\geq2.$$
\end{theorem}

\section*{Statements and Declarations}

\noindent\textbf{Competing interests} The author declares that he has no competing interests.


\begin{thebibliography}{99}

\bibitem{RACSAM6}
Lj. Aramba\v{s}i\'{c}, A. Guterman, B. Kuzma, R. Raji\'{c}, R. Tanaka, and S. Zhilina.
Geometric nonlinear classification of Hilbert $C^\ast$-modules based on strong Birkhoff-James orthogonality.
Rev. R. Acad. Cienc. Exactas F\'{i}s. Nat. Ser. A Mat. RACSAM, 118:151, 2024.
\url{https://doi.org/10.1007/s13398-024-01633-9}.

\bibitem{ArambasicBJMA2020}
Lj. Aramba\v{s}i\'{c}, A. Guterman, B. Kuzma, R. Raji\'{c}, and S. Zhilina.
Orthograph related to mutual strong Birkhoff-James orthogonality in $C^\ast$-algebras.
Banach J. Math. Anal., 14(4):1751--1772, 2020.
\url{https://doi.org/10.1007/s43037-020-00074-x}.

\bibitem{StrongMutualPreservers}
Lj. Aramba\v{s}i\'{c}, A. Guterman, B. Kuzma, R. Raji\'{c}, and S. Zhilina.
Operators preserving mutual strong Birkhoff-James orthogonality on $B(H)$.
Linear Algebra Appl., 624:27--43, 2021.
\url{https://doi.org/10.1016/j.laa.2021.04.003}.

\bibitem{Kuzma2021JMAA}
Lj. Aramba\v{s}i\'{c}, A. Guterman, B. Kuzma, R. Raji\'{c}, and S. Zhilina.
Symmetrized Birkhoff-James orthogonality in arbitrary normed spaces.
J. Math. Anal. Appl., 502(1):125203, 2021.
\url{https://doi.org/10.1016/j.jmaa.2021.125203}.

\bibitem{ArambasicSurvey}
Lj. Aramba\v{s}i\'{c}, A. Guterman, B. Kuzma, and S. Zhilina.
Birkhoff-James orthogonality: characterizations, preservers, and orthogonality graphs.
In R. M. Aron, M. S. Moslehian, I. M. Spitkovsky, and H. J. Woerdeman, editors,
Operator and Norm Inequalities and Related Topics, pages 255--302.
Birkh\"{a}user/Springer, Cham, 2022.
\url{https://doi.org/10.1007/978-3-031-02104-6_8}.

\bibitem{ArambasicAFA2014}
Lj. Aramba\v{s}i\'{c} and R. Raji\'{c}.
A strong version of the Birkhoff-James orthogonality in Hilbert $C^\ast$-modules.
Ann. Funct. Anal., 5(1):109--120, 2014.
\url{https://doi.org/10.15352/afa/1391614575}.

\bibitem{BhatiaSemrl}
R. Bhatia and P. \v{S}emrl.
Orthogonality of matrices and some distance problems.
Linear Algebra Appl., 287(1--3):77--85, 1999.
\url{https://doi.org/10.1016/S0024-3795(98)10134-9}.

\bibitem{Birkhoff}
G. Birkhoff.
Orthogonality in linear metric spaces.
Duke Math. J., 1(2):169--172, 1935.
\url{https://doi.org/10.1215/S0012-7094-35-00115-6}.

\bibitem{DolinarKuzma2016}
G. Dolinar and B. Kuzma.
Homomorphisms of commutativity relation.
Linear Multilinear Algebra, 64(5):897--922, 2016.
\url{https://doi.org/10.1080/03081087.2015.1066299}.

\bibitem{BJNorm}
A. Guterman, B. Kuzma, S. Singla, and S. Zhilina.
Birkhoff-James classification of norm's properties.
Adv. Oper. Theory, 9:43, 2024.
\url{https://doi.org/10.1007/s43036-024-00321-0}.

\bibitem{James1}
R. C. James.
Orthogonality in normed linear spaces.
Duke Math. J., 12(2):291--302, 1945.
\url{https://doi.org/10.1215/S0012-7094-45-01223-3}.

\bibitem{James2}
R. C. James.
Inner products in normed linear spaces.
Bull. Amer. Math. Soc., 53:559--566, 1947.
\url{https://doi.org/10.1090/S0002-9904-1947-08831-5}.

\bibitem{James3}
R. C. James.
Orthogonality and linear functionals in normed linear spaces.
Trans. Amer. Math. Soc., 61(2):265--292, 1947.
\url{https://doi.org/10.1090/S0002-9947-1947-0021241-4}.

\bibitem{Keckic2023JMAA}
D. J. Ke\v{c}ki\'{c} and S. Stefanovi\'{c}.
Isolated vertices and diameter of the BJ-orthograph in $C^\ast$-algebras.
J. Math. Anal. Appl., 528(1):127476, 2023.
\url{https://doi.org/10.1016/j.jmaa.2023.127476}.

\bibitem{KuzmaSushilPAMS}
B. Kuzma and S. Singla.
Birkhoff-James classification of finite-dimensional $C^\ast$-algebras.
Proc. Amer. Math. Soc., 153(4):1709--1721, 2025.
\url{https://doi.org/10.1090/proc/17153}.

\bibitem{KuzmaSushilFilomat}
B. Kuzma and S. Singla.
Non-linear classification of finite-dimensional simple $C^\ast$-algebras.
Filomat, 39(5):1465--1475, 2025.
\url{https://doi.org/10.2298/FIL2505465K}.

\bibitem{KuzmaSushilDebrecen}
B. Kuzma and S. Singla.
Classification of abelian finite-dimensional $C^\ast$-algebras by orthogonality.
Publ. Math. Debrecen, 108(1--2):45--73, 2026.
\url{https://doi.org/10.5486/PMD.2026.10162}.

\bibitem{KST}
B. Kuzma, S. Stefanovi\'{c}, and R. Tanaka.
On preservers of strong Birkhoff-James orthogonality between $C^\ast$-algebras.
arXiv:2510.24308v1 [math.OA] (2025).
\url{https://arxiv.org/abs/2510.24308v1}.

\bibitem{Magajna1993}
B. Magajna.
On the distance to finite-dimensional subspaces in operator algebras.
J. Lond. Math. Soc. (2), 47(3):516--532, 1993.
\url{https://doi.org/10.1112/jlms/s2-47.3.516}.

\bibitem{Stampfli}
J. G. Stampfli.
The norm of a derivation.
Pacific J. Math., 33(3):737--747, 1970.
\url{https://doi.org/10.2140/pjm.1970.33.737}.

\bibitem{Stefanovic}
S. Stefanovi\'{c}.
Diameter of the BJ-orthograph in finite-dimensional $C^\ast$-algebras.
Linear Multilinear Algebra, 73(4):718--728, 2025.
\url{https://doi.org/10.1080/03081087.2024.2382983}.

\end{thebibliography}
\end{document}